\documentclass[10pt]{amsart}
\usepackage{amssymb,amsmath,epsfig,graphics,latexsym,psfrag}
\usepackage[english]{babel}
\usepackage[all]{xy}
\newtheorem{theorem}{Theorem}[section]
\newtheorem{proposition}[theorem]{Proposition}
\newtheorem{corollary}[theorem]{Corollary}
\newtheorem{lemma}[theorem]{Lemma}
\newtheorem{claim}[theorem]{Claim}

\theoremstyle{definition}
\newtheorem{definition}[theorem]{Definition}

\theoremstyle{remark}
\newtheorem{remark}[theorem]{Remark}

\newtheorem{question}{Question}
\numberwithin{equation}{section}
\renewcommand{\t}{ \widetilde}
\renewcommand{\hat}{ \widehat}
\renewcommand{\b}{ \partial}
\newcommand{\Z}{\bf Z}
\newcommand{\R}{\bf R}

\newcommand{\Hi}{\bf H}

\renewcommand{\S}{\bf S}
\renewcommand{\l}{\langle}
\renewcommand{\r}{\rangle}
\newcommand{\e}{\varepsilon}

\renewcommand{\o}{\overline}

\newcommand{\co}{\colon\thinspace}
\renewcommand{\epsilon}{\varepsilon}
\renewcommand{\c}{\mathcal}

\usepackage{times}

\begin{document}
\sloppy

\title[]{bounded cohomology and $l_1$-homology   of  three-manifolds}

\author{Pierre Derbez}
\address{LATP, UMR 6632, 
Centre de Math\'ematiques et d'Informatique, 
Technopole de Chateau-Gombert, 
39, rue Fr\'ed\'eric Joliot-Curie -
 13453 Marseille Cedex 13}
\curraddr{}
\email{derbez@cmi.univ-mrs.fr}




\date{\today}

\begin{abstract}
 In this paper we define, for each aspherical orientable $3$-manifold $M$ endowed with a \emph{torus splitting} $\c{T}$, a $2$-dimensional fundamental $l_1$-class $[M]^{\c{T}}$ whose $l_1$-norm has similar properties as the Gromov simplicial volume of $M$ (additivity under torus splittings and isometry under finite covering maps).  Next, we use the Gromov simplicial volume of $M$ and the $l_1$-norm of $[M]^{\c{T}}$ to give a complete characterization of those nonzero degree maps $f\co M\to N$ which are homotopic to a ${\rm deg}(f)$-covering map. As an application we characterize those degree one maps $f\co M\to N$ which are homotopic to a homeomorphism in terms of bounded cohomology classes.
\end{abstract}
\maketitle

\vspace{-.5cm}

\section{Introduction}
 In order to characterize those nonzero degree maps $f\co M\to N$ between closed orientable aspherical $3$-manifolds which are homotopic to a finite covering we first need the simplicial volume of Gromov. 
 Recall that this invariant is obtained as follows: consider the $l_1$-norm on the real singular chains which induces, taking the infimum for all cycles representing a homology class, a $l_1$ semi-norm on the homology groups of a manifold. Then the simplicial volume of $M$ is the $l_1$ semi-norm of the fundamental class  $[M]$ corresponding to the orientation of $M$. 
This invariant is necessary since for any finite covering maps  $p\co\t{M}\to M$  then $(*)$ $$\|\t{M}\|=|{\rm deg}(p)|\|M\|$$ 
Equivalently, this equality means that $p$ induces an isometry $p_{\sharp}\co H_3(\t{M};{\R})\to H_3(M;{\R})$ with respect to the $l_1$ semi-norm. For hyperbolic $3$-manifolds equality $(*)$ is sufficient to characterize covering maps. More precisely, it follows from the Perelman geometrization of $3$-manifolds and from Gromov and Thurston's works  that if $M$ is a closed orientable hyperbolic $3$-manifold then any nonzero degree map $f\co M\to N$ into a closed orientable irreducible $3$-manifold is homotopic to a covering map if and only if $f_{\sharp}\co H_3(M;{\R})\to H_3(N;{\R})$ is an isometry.
However this condition is not sufficient to characterize covering maps between   non-hyperbolic $3$-manifolds. Roughly speaking this comes from the fact that the Gromov simplicial volume does not detect the "non-hyperbolic part" of an aspherical $3$-manifold.  
Thus one of the purposes of this paper is to construct for each aspherical oriented closed $3$-manifold $(M,\c{T}_M)$, where $\c{T}_M$ denotes the JSJ-family of canonical tori of $M$, a kind of \emph{secondary fundamental class} of $M$ denoted by $[M]^{\c{T}_M}$ which detects the non-hyperboic part of $M$. This class lies in the second $l_1$-homology group of $M$ denoted by   $H_2^{l_1}(M;{\R})$ and endowed with the $l_1$ semi-norm $\|.\|_1$. The  $l_1$-norm of $[M]^{\c{T}_M}$ together with the Gromov simplicial volume allow to characterize  finite covering maps. More precisely, the main result of this paper states as follows:

\begin{theorem}\label{localiso} Let  
 $f\co M\to N$ denote a nonzero degree map from an aspherical  closed oriented  $3$-manifold into an irreducible closed oriented $3$-manifold. If $M$ is not a $\t{{\rm SL}}(2,{\R})$-manifold then $f$    is homotopic to a ${\rm deg}(f)$-covering map iff  $\|f_{\sharp}([M])\|_1=\|[M]\|_1$ and $\|f_{\sharp}([M]^{\c{T}_M})\|_1=\|[M]^{\c{T}_M}\|_1$ where  $\c{T}_M$ denotes the JSJ-splitting  of $M$.
\end{theorem}

As a consequence of Theorem \ref{localiso} we deduce, using the duality between $l_1$-homology and bounded cohomology, a complete description of those degree one maps which are homotopic to a homeomorphism which answers to a question of M. Boileau. Denote by ${H}^2_b(M;{\R})$ the second bounded cohomology group of $M$ endowed with the semi-norm $\|.\|_{\infty}$ and by $\hat{H}^2_b(M;{\R})$, resp. $\hat{H}_2^{l_1}(M;{\R})$ the quotient space ${H}^2_b(M;{\R})/\ker\|.\|_{\infty}$, resp. ${H}_2^{l_1}(M;{\R})/\ker\|.\|_{1}$.
 
\begin{corollary}\label{rigidity}
A degree one map  $f\co M\to N$  between closed Haken manifolds  is homotopic to a homeomorphism iff $\|M\|=\|N\|$ and    $f$ induces an isometric isomorphism $f^{\sharp}\co(\hat{H}^2_b(N;{\R}),\|.\|_{\infty})\to(\hat{H}^2_b(M;{\R}),\|.\|_{\infty})$, resp. an isometry $f_{\sharp}\co(\hat{H}_2^{l_1}(M;{\R}),\|.\|_{1})\to(\hat{H}_2^{l_1}(N;{\R}),\|.\|_{1})$.
\end{corollary}

Recall that when $M$  admits a geometry ${\Hi}^3$ then M. Gromov and W.P. Thurston gave a characterization of local isometry using the Gromov simplicial volume.  When $M$ admits a geometry  $\t{\rm SL}(2,{\R})$ then Y. Rong solved this problem using the Seifert volume in \cite[Corollary 5.1]{Ro1} and when $M$  admits a geometry ${\R}^3$, ${\rm Nil}$ or ${\rm Sol}$ then S. Wang showed in \cite{W} that any nonzero degree map $f\co M\to N$ is homotopic to a covering map. On the other hand, when $M$ is a surface bundle over ${\S}^1$ (which covers the case where $M$ admits a ${\Hi}^2\times{\R}$-geometry) then M. Boileau and S. Wang gave a characterization of finite covering maps in terms of isometries with respect to the Thurston's norm (see \cite{BW}).

\section{Filling isometries and the secondary fundamental class}
In this section we define the objects and we state the intermediate results which will be used in the proof of Theorem \ref{localiso} and Corollary \ref{rigidity}. 

First of all recall recall that the $l_1$-homology groups of a topological space $X$ are obtained as follows. Denote by $C^{l_1}_{\ast}(X)$ the $l_1$-completion of the real singular chains $C_{\ast}(X)$. This means that 
$$C^{l_1}_{\ast}(X)=\left\{c=\sum_{i=1}^{\infty}a_i\sigma_i\ {\rm s.t.}\ \|c\|_1=\sum_{i=1}^{\infty}|a_i|<\infty\right\}$$
 where $a_i\in{\R}$ and $\sigma_i\co\Delta_{\ast}\to X$ is a singular simplex. We will denote by $S_{\ast}(X)$ the set of singular simplices. The topological dual of $C^{l_1}_{\ast}(X)$ is given by the set 
$$C^{\ast}_b(X)=\left\{w\in C^{\ast}(X)\ {\rm s.t.}\ \|w\|_{\infty}=\sup_{\sigma\in  S_{\ast}(X)}|\l w,\sigma\r|<\infty\right\}$$
Note that the $\b$ and $\delta$ operators are bounded so that $(C^{l_1}_{\ast}(X),\b)$ and $(C^{\ast}_b(X),\delta)$ are chain, resp. cochain, complexes.
 We denote by $H_{\ast}^{l_1}(X)$, resp. by $H^{\ast}_b(X)$, the homology, resp. cohomology, of this chain, resp. cochain, complex. The vector spaces $H_{\ast}^{l_1}(X)$ and $H^{\ast}_b(X)$ are endowed with the quotient semi-norm that we denote still by $\|.\|_1$ and $\|.\|_{\infty}$ respectively. Since a bounded operator has not necessarily a closed image then the above semi-norms are not norms in general.  Thus it will be convenient to consider the reduced $l_1$-homology and  bounded cohomology groups defined by $\hat{H}_{\ast}^{l_1}(X)=\ker\b/\o{{\rm Im}(\b)}=H_{\ast}^{l_1}(X)/\ker{\|.\|_1}$ and $\hat{H}^{\ast}_b(X)=\ker\delta/\o{{\rm Im}(\delta)}=H^{\ast}_b(X)/\ker{\|.\|_{\infty}}$. The evaluation map gives  Kronecker product 
$$\l .,.\r\co H^{\ast}_b(X)\otimes H_{\ast}^{l_1}(X)\to{\R}$$ 
which descends to a Kronecker product on reduced groups by the Holder inequality
$$\l .,.\r\co\hat{H}^{\ast}_b(X)\otimes\hat{H}_{\ast}^{l_1}(X)\to{\R}$$
\subsection{A filling isometry}
We begin with the construction of a filling  homomorphism. 
\begin{lemma}\label{rempli}
Let $(X,Y,A)$ be a triple of spaces where $Y\subset X$ (may be $X=Y$)  and each component of $A\subset Y$ has an amenable fundamental group. 

(i)  For any relative $2$-cycle $z$ in $(Y,A)$   there exists $u\in C_2^{l_1}(A)$ such that $z+u\in Z_2^{l_1}(Y)$.

(ii) Filling homomorphism. The map $z\mapsto z+u$ induces a filling homomorphism $$\Theta_Y\co H_2(Y,A)\to H_2^{l_1}(Y)\to H_2^{l_1}(X)$$ defined by $\Theta_Y([z])=[z+u]$.

(iii) Filling contraction. The filling homomorphism $\Theta_Y\co H_2(Y,A)\to H_2^{l_1}(X)$ satisfies $\|\Theta_Y\|\leq 1$ when $H_2(Y,A)$ is endowed with the $l_1$ semi-norm.

(iv) Naturality. Let $f\co(X,Y,A)\to(Z,W,B)$ be a continuous map of triple of spaces, where each component of $A$ and $B$ has amenable fundamental group. Then   the following diagram  is consistant.
$$\xymatrix{
H_2(Y,A;{\R}) \ar[r]^{\Theta_Y} \ar[d]_{f_{\sharp}} & H_2^{l_1}(X;{\R}) \ar[d]_{f_{\sharp}} \\
H_2(W,B;{\R}) \ar[r]^{\Theta_W}  & H_2^{l_1}(Z;{\R}) }
$$ 
\end{lemma}

In particular this result implies that $\Theta_Y$ induces a contraction homomorphism between reduced homology groups
$$\hat{\Theta}_Y\co\hat{H}_2(Y,A)=\frac{H_2(Y,A)}{\|.\|_1}\to\hat{H}_2^{l_1}(X)$$

Let $M$ denote a closed aspherical orientable $3$-manifold and let $\c{T}$ denote a  family of incompressible tori in $M$. We say that $\c{T}$ is a \emph{torus splitting} of $M$  if each component of  $M\setminus\c{T}$ is either a Seifert manifold or has a hyperbolic interior. From now on we assume that $M$ is endowed with a torus splitting $\c{T}$.   
 Let $(P,\b P)$ be a component of $M\setminus\c{T}$.  Let $\alpha\in H_2(P,\b P;{\R})$. Since each component of  $\b P$ has an amenable fundamental group then we can consider the filling homomorphism 
$$\hat{\Theta}_P\co\hat{H}_2(P,\b P;{\R})\to\hat{H}_2^{l_1}(M;{\R})$$
and the subspace of $H_2^{l_1}(M;{\R})$
$$\hat{H}_2^{l_1}\left(M^{\c{T}};{\R}\right)={\rm Vect}\left\l{\rm Im}(\hat{\Theta}_P)| P\ {\rm is\ a\ Seifert\ component\ of}\ M\setminus\c{T}\right\r$$  

Suppose that $P$ admits a Seifert fibration whose fiber is denoted by $h$. By a  \emph{horizontal surface} we mean a properly embedded  incompressible surface $F$ in $P$ which is transverse to $h$. Note that if $F$ is a horizontal surface in $P$ then the Seifert bundle $\eta$ induces an orbifold covering $\eta|F\co F\to\c{O}_P$ whose degree is denoted by $d_F\not=0$, where $\c{O}_P$ denotes the base $2$-orbifold. A horizontal surface will be termed \emph{minimal} if $|d_F|$ is minimal over all horizontal surfaces in $P$.
\begin{lemma}\label{generator}
Let $(P,h)$ be an aspherical oriented Seifert manifold with ${\S}^1$-fiber $h$ endowed with a fixed orientation $\mathfrak{o}(h)$. If $P$ is either Euclidean, ${\rm Nil}$ or a $\t{\rm SL}(2,{\R})$-manifold then $\hat{H}_2(P,\b P;{\R})=\{0\}$. If the interior of $P$ admits a ${\Hi}^2\times{\R}$ structure  then  $\hat{H}_2(P,\b P;{\R})\simeq{\R}$. Moreover, in the latter case,  if  $F$ and $F'$ are minimal surfaces in $P$ then $[F]=[F']$ in $\hat{H}_2(P,\b P;{\R})$ (where $[F]$, resp. $[F']$, corresponds to the orientation class of $F$, resp.  $F'$, so that $\mathfrak{o}(F)\times\mathfrak{o}(h)$, resp. $\mathfrak{o}(F')\times\mathfrak{o}(h)$, matches the given orientation of $P$)  and $\hat{H}_2(P,\b P;{\R})=\l[F]\r$. 
\end{lemma}
{\bf Notation.} For an oriented ${\Hi}^2\times{\R}$-manifold $P$ we will denote by $\alpha_P$ the class of any minimal surface in $\hat{H}_2(P,\b P;{\R})$ with the  convention for orientations of Lemma \ref{generator}. For Euclidean, ${\rm Nil}$ or a $\t{\rm SL}(2,{\R})$-manifold then we set $\alpha_P=0$.
We denote by $\hat{H}^+_2(P,\b P;{\R})$ the set defined by $$\{\alpha\in\hat{H}_2(P,\b P;{\R}) |  \ \alpha=\xi_.\alpha_{P},\ \xi\geq 0\}$$ When $M$ denotes an oriented aspherical $3$-manifold endowed with a torus splitting $\c{T}$ each Seifert component $P$ of $M\setminus\c{T}$ is oriented by $M$ and we  fix an orientation  $\mathfrak{o}(h)$ for the fiber $h$ of each $P$. In this case we say that $M$ is \emph{framed}. We denote  by $\hat{H}^{l_1,+}_2(M^{\c{T}};{\R})$ the set defined by  $$\{\alpha\in\hat{H}_2^{l_1}(M^{\c{T}};{\R}) |  \ \alpha=\sum\xi_i.\hat\Theta_{P_i}(\alpha_{P_i}),\ \xi_i\geq 0\}$$ where the $P_i$'s run over the Seifert components of $M\setminus\c{T}$. Notice that since the only aspherical Seifert fibered spaces admitting at least two non-isotopic  fibrations are the Euclidean manifolds then the above notions are well-defined by Lemma \ref{generator}.
Some geometric properties of the map $\hat{\Theta}_P$ are reflected in the following
\begin{theorem}\label{1}
Let $(M,\c{T})$ be a closed aspherical orientable $3$-manifold endowed with a torus splitting  and denote by  $P_1,...,P_l$ the Seifert components of $M\setminus\c{T}$.  Then 

(i) Isometry: the filling homomorphism $\hat{\Theta}_{P_i}\co\hat{H}_2(P_i,\b P_i;{\R})\to \hat{H}_2^{l_1}(M;{\R})$ is an isometry with respect to the $l_1$-norms. Moreover for any $\alpha\in\hat{H}_2(P_i,\b P_i;{\R})$ then $$\|\alpha\|_1=\|\hat\Theta_{P_i}(\alpha)\|_1=|\xi_{\alpha}|\|F_i\|$$ where $\xi_{\alpha}$ is the real number such that $\alpha=\xi_{\alpha}.\alpha_{P_i}$ and $F_i$ is a minimal surface in $P_i$. In particular, if $\alpha$ is represented by an  incompressible connected surface  $\c{F}$ then $$\|\alpha\|_1=\|\hat\Theta_{P_i}(\alpha)\|_1=\|\c{F}\|$$

(ii) Additivity: for any $l$-uple $(\alpha_1,...,\alpha_l)\in\hat{H}^+_2(P_1,\b P_1;{\R})\times...\times\hat{H}^+_2(P_l,\b P_l;{\R})$ we have 
$$\|\hat\Theta_{P_1}(\alpha_1)+...+ \hat\Theta_{P_l}(\alpha_l)\|_1=\|\hat\Theta_{P_1}(\alpha_1)\|_1+...+ \|\hat\Theta_{P_l}(\alpha_l)\|_1$$
\end{theorem}
We end this section with a result which describes the metric behavior of finite covering maps. To this purpose note that throughout this paper we adopt the following convention for the orientations. Let $(M,\c{T})$ be a closed aspherical framed $3$-manifold endowed with a torus splitting and denote by $p\co(\t{M},\t{\c{T}})\to(M,\c{T})$ a finite covering endowed with a torus splitting defined by $\t{\c{T}}=p^{-1}(\c{T})$. Let $(\Sigma,h)$ denote a Seifert piece of $M\setminus\c{T}$, where $h$ denotes the fiber of $\Sigma$ and let $(\Sigma_1,h_1),...,(\Sigma_l,h_l)$ the components of $p^{-1}(\Sigma)$ where $h_i$ denotes the fiber of $\Sigma_i$  so that $p|\Sigma_i\co\Sigma_i\to\Sigma$ is a fiber preserving map (such a Seifert fibration always exists on $\Sigma_i$ by \cite{JS}).   Then we orient the fibers $h_i$ so that $p|h_i\co h_i\to h$ is orientation preserving. In this case we say that $\t{M}$ is endowed with the framing induced by that of $M$.

\begin{proposition}\label{2} Let $(M,\c{T})$ be a closed aspherical framed $3$-manifold endowed with a torus splitting.
 Any finite covering  $p\co\t{M}\to M$    induces   isometries $p_{\sharp}\co H_3\left(\t{M};{\R}\right)\to H_3\left(M;{\R}\right)$ and    $$p_{\sharp}|\hat{H}_2^{l_1,+}\left(\t{M}^{\t{\c{T}}};{\R}\right)\co\hat{H}_2^{l_1,+}\left(\t{M}^{\t{\c{T}}};{\R}\right)\to\hat{H}_2^{l_1}\left(M;{\R}\right)$$ where $\t{\c{T}}$ is the torus splitting of $\t{M}$ equal to $p^{-1}(\c{T})$ and where $\t{M}$ is endowed with the framing induced by that of $M$.  
\end{proposition}
Note that the above \emph{covering property} gives rise to the following 
\begin{question}\label{iff}
Let $M$ be a closed orientable aspherical $3$-manifold and let $p\co\t{M}\to M$ denote a Haken finite covering of $M$.

(i) Does the covering induce an isometry $p_{\sharp}|\hat{H}_2^{l_1,+}\left(\t{M}^{\t{\c{T}}};{\R}\right)\co\hat{H}_2^{l_1,+}\left(\t{M}^{\t{\c{T}}};{\R}\right)\to\hat{H}_2^{l_1}\left(M;{\R}\right)$ for any torus splitting $\t{\c{T}}$  of $\t{M}$? 

(ii) More generally does the covering induce an isometry $p_{\sharp}\co\hat{H}_2^{l_1}\left(\t{M};{\R}\right)\to\hat{H}_2^{l_1}\left(M;{\R}\right)$?
\end{question}
\subsection{The secondary fundamental class of a manifold}
We   can now define the secondary fundamental class   of a closed orientable aspherical $3$-manifold $M$. Consider  a torus splitting $\c{T}$ of $M$ and denote by $P_1,...,P_k$ the components of $M\setminus\c{T}$ supporting a Seifert fibration. Then we set $$\|M\|^{\c{T}}=\|\hat{\Theta}_{P_1}(\alpha_{P_1})+...+\hat{\Theta}_{P_k}(\alpha_{P_k})\|_1$$
The class $\hat{\Theta}_{P_1}(\alpha_{P_1})+...+\hat{\Theta}_{P_k}(\alpha_{P_k})\in\hat{H}_2^{l_1,+}\left({M}^{{\c{T}}};{\R}\right)$ will be denoted by $[M]^{\c{T}}$ and will be termed the secondary fundamental class of $(M,\c{T})$.

\begin{question}
Is it possible to compare $[M]^{\c{T}_1}$ and $[M]^{\c{T}_2}$ in $\hat{H}_2^{l_1}\left({M};{\R}\right)$ when $\c{T}_1$ and $\c{T}_2$ are two distinct torus splittings of $M$? 
\end{question}

\subsection{Organization of the paper}
In section 3 we prove Lemmas \ref{rempli} and \ref{generator} and we recall some results related to the duality between bounded cohomology and $l_1$-homology. Section 4 will be devoted the construction of bounded $2$-cocyles which roughly speaking measure the area of the horizontal surfaces passing through the Seifert pieces of a manifold.  This kind of cocyles will be used the estimate  the $l_1$-norm certain $l_1$-homology classes (see Proposition \ref{bounded}). In section 5 we prove Theorem \ref{1} and Proposition \ref{2} and section 6 is devoted to the proof of Theorem \ref{localiso} and Corollary \ref{rigidity}.

\section{Filling homomorphisms and Duality between bounded cohomology and $l_1$-homology}
\begin{proof}[Proof of Lemma \ref{rempli}]
Let $\alpha\in H_2(Y,A)$ and choose  a relative $2$-cycle $z$ in $(Y,A)$ which represents $\alpha$. Then $\b z\in Z_1(A)$ and since $H^{l_1}_1(A)=0$ then there exists $u\in C_2^{l_1}(A)$ such that $\b u=-\b z$ in such a way that $z+u\in Z_2^{l_1}(Y)$. This proves point (i).

On the other hand the class of $z+u$ does not depend on the choice of  $u$. Indeed let $u_1,u_2$ be elements of $C_2^{l_1}(A)$ such that $\b u_1=\b u_2=-\b z$. Then $(z+u_1)-(z+u_2)=u_1-u_2\in Z_2^{l_1}(A)$. Then there exists $w\in C^{l_1}_3(A)$ such that $\b w=(z+u_1)-(z+u_2)$. 

Moreover the class of $z+u$ does not depend on the choice of  the representant $z$. Indeed, let $z$ and $z'$ be two relative $2$-cycles  which represent $\alpha$. Then there exists $p\in C_2(A)$ and $v\in C_3(Y)$ such that $z-z'=\b v+p$. Let $u$ and $u'$ in $C^{l_1}_2(A)$ such that $\b u=-\b z$ and $\b u'=-\b z'$. Then $(z+u)-(z'+u')=\b v+(p+u-u')$. Thus $p+u-u'\in Z^{l_1}_2(A)$. Then there exists $w\in C^{l_1}_3(A)$ such that $\b w=p+u-u'$. Thus the map $z\mapsto z+u$ gives rise to a  homomorphism
$$\Theta_Y\co H_2(Y,A)\to H_2^{l_1}(Y)\to H_2^{l_1}(X)$$ defined by $\Theta_Y([z])=[z+u]$. This proves point (ii). 

Let $\alpha\in H_2(Y,A)$ and fix   $\e>0$. Then by \cite[Equivalence Theorem]{G} there exists a representative $z$ of $\alpha$ such that $\|z\|_1\leq\|\alpha\|_1+\e$ and $\|\b z\|_1\leq\e$. By the Uniform Boundary Condition (see \cite[Theorem 2.3]{MM}) there exists a constant $K>0$ which only depends on the dimension and $u\in C^{l_1}_2(A)$ such that $\b u=-\b z$ and $\|u\|_1\leq K\|\b z\|_1<K\e$. This implies, passing to the limits, that $\|\Theta_Y(\alpha)\|_1\leq\|\alpha\|_1$. This proves point (iii).

Let $\alpha$ be an element of $H_2(Y,A;{\R})$. Then $f_{\sharp}(\Theta_Y\alpha)=[f_{\sharp}(z)+f_{\sharp}(u)]$ where $z$ is a relative $2$-cycle representing $\alpha$ and $u$ is a $l_1$ $2$-chain in $A$ such that $\b u=-\b z$. Since  $\b f_{\sharp}u=-\b f_{\sharp}z$ and since $f_{\sharp}u$ is a $l_1$ $2$-chain in $B$ then $[f_{\sharp}(z)+f_{\sharp}(u)]=\Theta_W(f_{\sharp}(\alpha))$. This proves of point (iv).

\end{proof}
\begin{proof}[Proof of Lemma \ref{generator}]

Let $(P,h)$ be a Seifert oriented $3$-manifold with a fixed fibration $h$. Let $p\co\t{P}\to P$ denote a finite covering of $P$. Then $p$ induces an epimorphism ${H}_2(\t{P},\b\t{P})\to{H}_2(P,\b P)$ and thus passing to the quotient we get an epimorphism $\hat{H}_2(\t{P},\b\t{P})\to\hat{H}_2(P,\b P)$.

If $P$ admit a geometry ${\rm Nil}$, ${\rm R}^3$ or $\t{\rm SL}(2,{\R})$ then it admits a finite covering $\t{P}$ which is either a torus bundle (over ${\S}^1$ or $I$) or a circle bundle over a hyperbolic surface with non-zero Euler number.   In any case $\hat{H}_2(\t{P},\b\t{P})=\{0\}$. Indeed in these cases $\t{P}$ contains no incompressible surfaces with negative Euler characteristic. This proves that   $\hat{H}_2({P},\b{P})$ is trivial.

 Hence from now on we assume that ${\rm int}(P)$ admits a ${\Hi}^2\times{\R}$-geometry.  Let $p\co\t{P}\to P$ denotes a finite regular covering of $P$ homeomorphic to a product $\t{F}\times{\S}^1$.   Note that since $\t{P}$ has an orientable base then by \cite[Lemma 6]{WZ} ${H}_2(\t{P},\b\t{P})$ is generated by any minimal surface together with vertical surfaces (which are either tori or vertical annuli) and thus  $\hat{H}_2(\t{P},\b\t{P})$ is generated by any minimal surface. This proves  that $\hat{H}_2(\t{P},\b\t{P})={\R}$ and thus $\hat{H}_2(P,\b P)={\R}$.

Let us check the second statement of the lemma. Let $F_0$ and $F_1$ denote two minimal surfaces in $P$ and denote by $\t{F}_0$, resp. $\t{F}_1$ the spaces $p^{-1}(F_0)$, resp. $p^{-1}(F_1)$. Since $F_0$ and $F_1$ are minimal then in particular $\chi(F_0)=\chi(F_1)$ and thus $\chi(\t{F}_0)=\chi(\t{F}_1)$ which implies, since $\t{P}$ is a product that $[\t{F}_0,\b\t{F}_0]=[\t{F}_1,\b\t{F}_1]$ in $\hat{H}_2(\t{P},\b\t{P})$. Since for each $i=0,1$ we have $p_{\sharp}([\t{F}_i,\b\t{F}_i])={\rm deg}(p)[F_i,\b F_i]$ then $[F_0,\b F_0]=[F_1,\b F_1]$ which proves the lemma.

\end{proof}

 Throughout this paper we will need the following general results which come from the duality between $l_1$-chains and bounded cochains.
\begin{lemma}\label{duality}
The Kronecker product between $l_1$-homology and bounded cohomology gives rise to a surjective bounded operator
$$\Phi\co\hat{H}^{\ast}_b(X)\to(\hat{H}_{\ast}^{l_1}(X))'$$ where $(\hat{H}_{\ast}^{l_1}(X))'$ denotes the space of continuous linear forms on $\hat{H}_{\ast}^{l_1}(X)$, such that $\|\Phi\|=1$ and  for any $\varphi\in(\hat{H}_{\ast}^{l_1}(X))'$ there exists $\beta\in\hat{H}^{\ast}_b(X)$ such that $\Phi(\beta)=\varphi$ and $\|\beta\|_{\infty}=\|\varphi\|_{\infty}$.
\end{lemma}
\begin{proof}
This follows directly from the Holder inequality combined with the Hahn-Banach Theorem.
\end{proof}
Let $\Gamma$ be a group acting by homeomorphisms on a topological space $X$. We denote by $\Gamma\hat{H}^{l_1}_{\ast}(X)=\{\alpha\in\hat{H}^{l_1}_{\ast}(X)\ {\rm s.t.} \ g_{\sharp}(\alpha)=\alpha\ {\rm when} \ g\in\Gamma\}$.
\begin{lemma}\label{galois}
Let $p\co\t{X}\to X$ be a regular covering map with finite Galois group $\Gamma$. Then the induced homomorphism $p_{\sharp}\co\Gamma\hat{H}^{l_1}_{\ast}(\t{X})\to\hat{H}^{l_1}_{\ast}(X)$ is an  isometry.
\end{lemma}
\begin{proof}
Consider the averaging retraction $A\co C_b^{\ast}(\t{X})\to C_b^{\ast}(X)$ defined by  
$$\l A(\gamma),\sigma\r=\frac{\sum_{g\in\Gamma}\l g^{\ast}\gamma,\t{\sigma}\r}{{\rm Card}(\Gamma)}$$ where $\t{\sigma}\co\Delta^{\ast}\to\t{X}$ denotes a lifting of $\sigma\co\Delta^{\ast}\to X$. This definition does not depend one the choice of the lifting $\t{\sigma}$ since the covering is regular. By construction $A$ satisfies the identity $A\circ p^{\ast}={\rm Id}$ and  commutes with the differentials so that it induces a homomorphism $\hat{A}\co H_b^{\ast}(\t{X})\to H_b^{\ast}(X)$. Then at the $H_b^{\ast}$-level we still have the identity $\hat{A}\circ p^{\ast}={\rm id}$ and $\|\hat{A}\|\leq 1$. Notice that $\hat{A}$ induces a homomorphism $\hat{H}_b^{\ast}(\t{X})\to\hat{H}_b^{\ast}(X)$ still denoted by $\hat{A}$.

Let $\alpha=\o{[z]}\in\Gamma\hat{H}^{l_1}_{\ast}(\t{X})$, where $-\co{H}^{l_1}_{\ast}(\t{X})\to\hat{H}^{l_1}_{\ast}(\t{X})$ denotes the natural quotient homomorphism and $z$ is a $l_1$-cycle. If $\alpha=\o{[z]}\not=0$ then by Lemma \ref{duality}, there exists $\beta=\o{[\gamma]}\in \hat{H}^{\ast}_b(\t{X})$ such that $\l\beta,\alpha\r=1$ and $\|\beta\|_{\infty}=\frac{1}{\|\alpha\|_1}$. 
Since $\alpha$ is $\Gamma$-invariant then by the definition of the averaging we have the equalities 
$$\left\l\hat{A}(\beta),p_{\sharp}(\alpha)\right\r=\left\l A(\gamma),p_{\sharp}(z)\right\r=\frac{1}{{\rm Card}(\Gamma)}\sum_{g\in\Gamma}\left\l g^{\ast}(\gamma),z\right\r=\left\l\gamma,z\right\r=\left\l\beta,\alpha\right\r=1$$
 and thus 
$$\|p_{\sharp}(\alpha)\|_1\geq\frac{1}{\|\hat{A}(\beta)\|_{\infty}}\geq\frac{1}{\|\beta\|_{\infty}}=\|\alpha\|_1$$ Since the inequality $\|p_{\sharp}\|\leq 1$ is always true this proves the lemma. 
\end{proof}

\section{Bounded $2$-cocyles measuring the horizontal areas}
 Let  $M$ be an  orientable closed aspherical $3$-manifold endowed with a torus splitting $\c{T}$. A Seifert piece $P$ of $(M,\c{T})$ will be termed a \emph{product component} if it is homeomorphic to a product $F\times{\S}^1$, where $F$ is a surface  whose interior admits a  hyperbolic structure.   The purpose of this section is to construct for each "product component" $P$ of $M\setminus\c{T}$,  a bounded $2$-cocyle in $M$ which measures the hyperbolic area of horizontal surfaces of $P$. More precisely the main result of this section states as follows:

 \begin{proposition}\label{bounded}
 For each product component $P=F\times{\S}^1$ of $M\setminus\c{T}$ there exists a non-trivial bounded $2$-cocyle $\Omega_P$ in $M$ satisfying the following properties:

(i) $i^{\ast}(\Omega_P)$ is a relative $2$-cocycle in $(P,\b P)$ where $i\co P\hookrightarrow M$ denotes the natural inclusion,

(ii) for any connected horizontal surface $\c{F}$  in $P$ then $$\left|\left\l[\Omega_P],\hat{\Theta}_P([\c{F}])\right\r\right|=\left|\l i^{\ast}(\Omega_P),z_{\c{F}}\r\right|={\rm Area}(\c{F})$$ where $z_{\c{F}}$ is a relative $2$-cycle representing the fundamental class of $\c{F}$  and ${\rm Area}(\c{F})$ denotes the area of ${\rm int}(\c{F})$ endowed with a complete hyperbolic metric. 

 Let $(P_i)_{i\in I}$ be a family of pairwise distinct product components of $M\setminus\c{T}$  and for each $i\in I$ denote by $k_i\co P_i\to M$ the canonical inclusion. Then 
 
 (iii) $k^{\ast}_i(\Omega_{P_j})=0$ if $i\not=j$ and 
 $$\left\|\sum_{i\in I}[\Omega_{P_i}]\right\|_{\infty}=\pi$$
where $[\Omega_{P_i}]$ denotes the class of $\Omega_{P_i}$ in $\hat{H}^2_b(M;{\R})$ for each $i\in I$. 
 \end{proposition}

 To prove this proposition, we first need to construct a chain map: the \emph{reduction}, used in \cite{FS} (see section \ref{red}). Next we \emph{straight horizontally} the reduced chains which meet essentially $P$ to define a bounded $2$-cochain $\Omega_P$ (see section \ref{horizontal}). The proof of Proposition \ref{bounded} will occupy section \ref{proofbounded}. From now on, we denote by $\Delta^n$ the standard $n$-simplex defined by 
 $$\Delta^n=[v_0,...,v_n]=\left\{\sum_{i=0}^{n}t_iv_i, \sum_{i=0}^{n}t_i=1, t_i\geq 0\right\}$$ where $v_i=(0,...,1,...0)\in{\R}^{n+1}$. We denote by $V(\Delta^n)=\{v_0,...,v_n\}$ the vertices of $\Delta^n$ and by $\Delta^{n-1}_i$ the $i$-th $(n-1)$-face of $\Delta^n$ defined by $\Delta_i^{n-1}=[v_0,...,\hat{v}_i,...,v_n]$. Note that $\Delta^n$ is oriented by the order of its vertices in such a way that $(v_1-v_0,...,v_n-v_0)$ is a positive frame.  
\subsection{Reduction operator in singular homology}\label{red}
Let $M$ be a closed aspherical orientable $3$-manifold endowed with an amenable splitting $\c{T}$. Denote by $P_1,...,P_l$ the  components of $M\setminus\c{T}$. As in \cite{FS}, we consider a chain map $\rho\co C_{n}(M)\to C_{n}(M)$ defined as follows:

If $n=0$ then we set $\rho=\mathbf{1}$.

If $n=1$ let $\tau\co[v_0,v_1]\to M$ be a 1-simplex. Since $\c{T}$ is incompressible, the map $\tau$ is homotopic rel. $\{v_0,v_1\}$ to a \emph{reduced} $1$-simplex i.e. a map $\tau_1\co[v_0,v_1]\to M$ such that either
$\tau_1([v_0,v_1])\subset\c{T}$ or 
 $\tau_1|(v_0,v_1)$ is transverse to $\c{T}$ and for each component $J$ of $\tau_1^{-1}(P_i)$ then $\tau_1|J$ is not homotopic rel. $\b J$ into $\b P_i$.
Then we set $\rho(\tau)=\tau_1$ and we extend $\rho$   by linearity. 

If $n=2$ let $\sigma\co\Delta^2=[v_0,v_1,v_2]\to M$ be a $2$-simplex.  Then $\sigma$ is homotopic rel. $V(\Delta^2)=\{v_0,v_1,v_2\}$ to a reduced $2$-simplex $\sigma_1$ such that either 
 $\sigma_1(\Delta^2)\subset\c{T}$ or 
 $\sigma_1|{\rm int}(\Delta^2)$ is transverse to 
$\c{T}$, the $1$-simplex $\sigma_1|e$ is reduced for each edge $e$ of $\Delta^2$ and  $\sigma_1^{-1}(\c{T})$ contains no loop components. Thus if $J$ is a component of $\sigma_1^{-1}(\c{T})$ such that $J\cap{\rm int}(\Delta^2)\not=\emptyset$ then $J$ is a proper arc in $\Delta^2$ connecting two distinct edges.   
Then we set $\rho(\sigma)=\sigma_1$ and we extend $\rho$   by linearity. 

\begin{remark}\label{ideal}
Let $\sigma\co\Delta^2\to M$ be a reduced $2$-simplex.  Denote by $D$ a component of $\sigma^{-1}({\rm int}(P))$, where $P$ is a component of $M\setminus\c{T}$. Each component  of $V(\Delta^2)\cap D$ will be termed a vertex of $D$ and each component of $\sigma^{-1}(\b P)\cap\o{D}$ will be termed an \emph{ideal vertex} of $D$. Thus an ideal vertex of $D$ is either a vertex or an edge or a proper arc of $\Delta^2$ connecting two distinct edges. Let us denote by $V_{\infty}(D)$ the set of ideal vertices of $D$ and by $V(D)$ the set of all  vertices and ideal vertices of $D$. Note that if $D\not=\emptyset$ then $2\leq{\rm Card}(V(D))\leq 3$.   
\end{remark}
\begin{remark}\label{core}
Suppose that $\sigma\co\Delta^2\to M$ is a reduced $2$-simplex.  If   $\sigma(e)$ is not contained in $\c{T}$ for any edge $e$ of $\Delta^2$ then there exists a unique component, denoted by ${\rm Core}(\sigma)$, of $\Delta^2\setminus\sigma^{-1}(\c{T})$ which meets the three edges of $\Delta^2$ (see \cite{FS}). On the other hand,  if $e_1$ and $e_2$ denote two distinct edges of $\Delta^2$ such that $\sigma(e_1)\subset\c{T}$ and $\sigma(e_2)\subset\c{T}$ then $\sigma(\Delta^2)\subset\c{T}$ (this follows directly from the reduction hypothesis of $\sigma$).  
\end{remark}

If $n=3$ let $\sigma\co\Delta^3=[v_0,v_1,v_2,v_3]\to M$ be a $3$-simplex.  Then $\sigma$ is homotopic rel. $V(\Delta^3)=\{v_0,v_1,v_2,v_3\}$ to a reduced $3$-simplex $\sigma_1$ such that either 
$\sigma_1(\Delta^3)\subset\c{T}$ or 
$\sigma_1|{\rm int}(\Delta^3)$ is transverse to 
$\c{T}$, the $2$-simplex $\sigma_1|\Delta^2_i$ is reduced for each face $\Delta^2_i$ of $\Delta^3$  and  if $D$ is a component of $\sigma_1^{-1}(\c{T})$ such that $D\cap{\rm int}(\Delta^3)\not=\emptyset$ then $D$ is either a  normal triangle or  a normal rectangle (see figure 1). 

Then we set $\rho(\sigma)=\sigma_1$ and we extend $\rho$   by linearity. For $n=0,1,2,3$ we denote by $C_{n}^{\rm red}(M)$ the image of $C_{n}(M)$ under $\rho$. Notice that the reduction is stable under the $\b$-operator. 

\begin{remark}\label{frontier}
Let $\sigma$ denote a reduced $3$-simplex and let $\nabla$ be a component of $\sigma^{-1}({\rm int}(P))$, where  $P$ is a component of $M\setminus\c{T}$. Each component  of $V(\Delta^3)\cap \nabla$ will be termed a vertex of $\nabla$ and each component of $\sigma^{-1}(\b P)\cap\o{\nabla}$ which is a normal triangle or a vertex of $\Delta^3$ will be termed an \emph{ideal vertex} of $\nabla$. Denote by $V_{\infty}(\nabla)$ the ideal vertices of $\nabla$ and by $V(\nabla)$ all the  vertices and ideal vertices of $\nabla$.  The set $\sigma^{-1}(\b P)\cap\o{\nabla}$ consists of either 

(i) exactly one $2$-face $\Delta^2_i$ of $\Delta^3$ and at most one ideal vertex of $\nabla^3$, or

(ii)  exactly two edges of $\Delta^3$ with no common vertices, or

(iii) exactly one edge and at most one normal rectangle or at most two ideal vertices, or

(iv) no faces and no edges components and   at most four ideal vertices or two rectangles or one rectangle and at most two ideal vertices. 

\end{remark}

\begin{center}
\psfrag{vi}{$v_i$} \psfrag{vj}{$v_j$} \psfrag{vk}{$v_k$} \psfrag{vl}{$v_l$}
\psfrag{wi}{$w_i$} \psfrag{wj}{$w_j$} \psfrag{wk}{$w_k$} \psfrag{wl}{$w_l$}
\psfrag{R}{$R$} 
\includegraphics{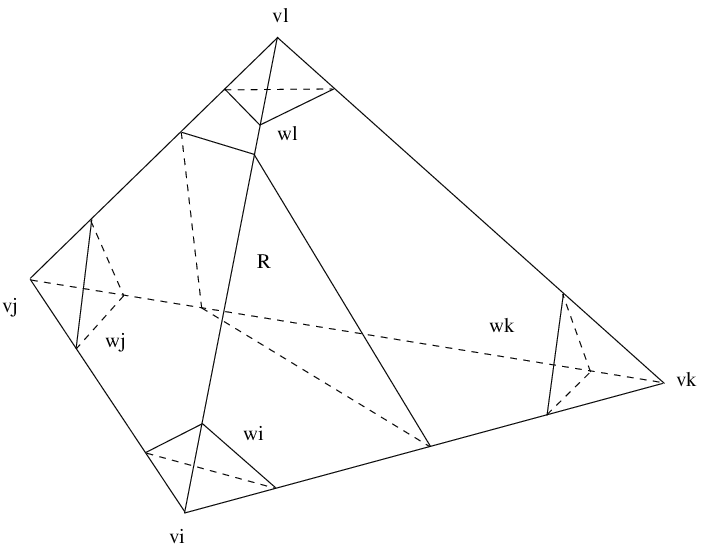}
\centerline{Figure 1: Normal triangles and rectangles}
\end{center}
\vskip 0.5 truecm

\subsection{Horizontal area}\label{horizontal} Suppose now that $M$ is endowed with a torus splitting $\c{T}$. Let $P$ be a product component of $M\setminus\c{T}$   identified with a product  $F\times{\S^1}$, where $F$ is an orientable  surface  whose interior $W$ admits a hyperbolic structure. Denote by $q\co P\to F$ the projection to the first factor.

Let $\sigma\co\Delta^2\to M$ be a reduced $2$-simplex and   let $D$ be a component of $\sigma^{-1}({\rm int}(P))$. Denote by $\sigma_D$ the map $\sigma|D\co D\to{\rm int}(P)$ and  by $\t{q\sigma_D}\co D\to{\Hi}^2$ a lifting of $q\sigma_D\co D\to W$ into the universal covering of $W$.   It follows from Remark \ref{ideal} that $D$ has $2$ or $3$ vertices (including ideal vertices). Denote by $V(D)=\{w_0,w_1,w_2\}$ the set of vertices of $D$ (may be $w_i=w_j$ for some $i\not=j$) in such a way that if $D={\rm Core}(\sigma)$ then each $w_i$ corresponds to $v_i$ (see figure 2).

\begin{center}
\psfrag{v0}{$v_0$} \psfrag{v1}{$v_1$} \psfrag{v2}{$v_2$}
\psfrag{w0}{$w_0$} \psfrag{w1}{$w_1$} \psfrag{w2}{$w_2$} \psfrag{w3}{$w_3$}
\psfrag{w0w1}{$w_0=w_1$} \psfrag{v1w1}{$v_1=w_1$}
\psfrag{I}{$(I)$} \psfrag{II}{$(II)$} \psfrag{III}{$(III)$}
\psfrag{IV}{$(IV)$} \psfrag{V}{$(V)$}
\psfrag{C}{${\rm Core}(\sigma)$} \psfrag{cvide}{${\rm Core}(\sigma)=\emptyset$}
\psfrag{D}{$D$}
\includegraphics{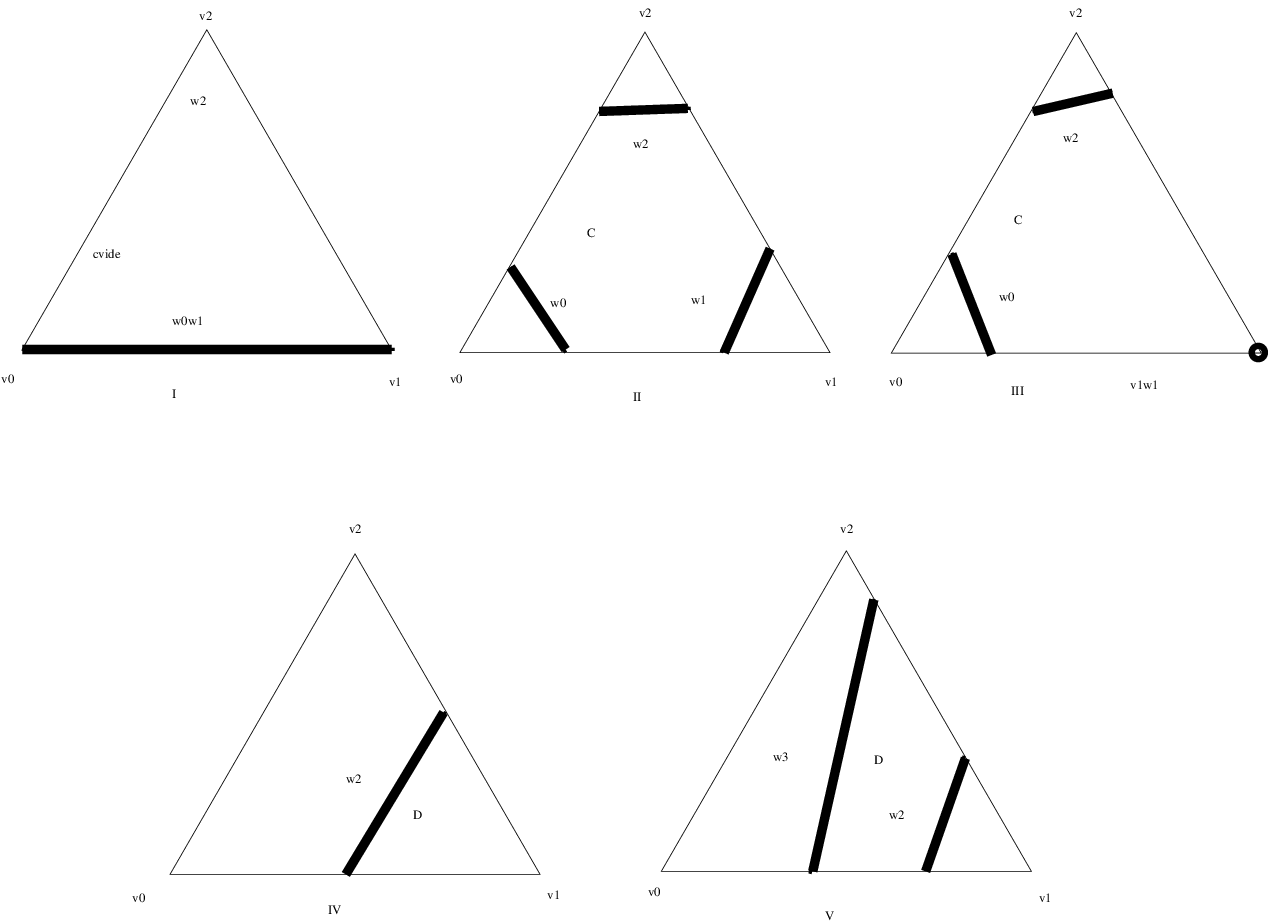}
\centerline{Figure 2}
\end{center}
\vskip 0.5 truecm

For each $i=0,1,2$, the vertex $w_i$ has a corresponding image $w_i^{\infty}$ in $\o{\Hi}^2={\Hi}^2\cup\b_{\infty}{\Hi}^2$ defined as follows:

If $w_i\not\in V_{\infty}(D)$ then $w_i^{\infty}=\t{q\sigma_D}(w_i)$. If $w_i\in V_{\infty}(D)$ denote by $w_j$ a vertex of $D$ distinct from $w_i$. Choose a point $x_i$ in $w_i$ and $x_j$ in $w_j$ and consider the open geodesic segment $(x_i,x_j)$ in $\Delta^2$.   Then there exists $\e>0$ and a horodisk $C_i$ in ${\Hi}^2$ such that $\t{q\sigma_D}((x_i,x_i+\e(x_j-x_i)])\in C_i$ and in this case  $w_i^{\infty}$ denotes the center of $C_i$ (which is defined as the contact point of $C_i$ with $\b_{\infty}{\Hi}^2$, see figure 3). Note that $w_i^{\infty}$ does not depend on the choice of the points $x_i$ and $x_j$ and does not depend on the choice of $w_j\not=w_i$ since in any case it follows from our construction and from the reduction hypothesis on $\sigma$ that 
$$\lim_{t\to 1^-}\t{q\sigma_D}(tx_i+(1-t)x_j)=w_i^{\infty}$$ for any $x_i\in w_i$, $x_j\in w_j$ and $w_j\not=w_i$.
On the other hand, since $\sigma$ is reduced, if $w_i\not=w_j$ where $w_i$ or $w_j$ is an ideal vertex then $w_i^{\infty}\not=w_j^{\infty}$. The set $\{w_0^{\infty},w_1^{\infty},w_2^{\infty}\}$ will be termed the $\sigma$-image of the (ideal) vertices of $D$.

\begin{center}
\psfrag{v0}{$v_0$} \psfrag{v2w2}{$v_2=w_2$} \psfrag{v1}{$v_1$} \psfrag{v2}{$v_2$}
\psfrag{w0}{$w_0$} \psfrag{w1}{$w_1$} \psfrag{w2}{$w_2$} \psfrag{w3}{$w_3$}
\psfrag{w0w1}{$w_0=w_1$} \psfrag{q}{$\t{q\sigma_D}$}
\psfrag{v0v2}{$\t{q\sigma_D}([v_0,v_2])$}  \psfrag{v1v2}{$\t{q\sigma_D}([v_1,v_2])$}
\psfrag{winf}{$w_2^{\infty}$}
\psfrag{c}{$C$}
\psfrag{infini}{$w_0^{\infty}=w_1^{\infty}$}
\psfrag{H}{${\Hi}^2$}
\psfrag{I}{$(I)$} \psfrag{II}{$(II)$} \psfrag{III}{$(III)$}
\psfrag{x}{$w_0^{\infty}$} \psfrag{y}{$w_1^{\infty}$} \psfrag{z}{$w_2^{\infty}$}
\includegraphics{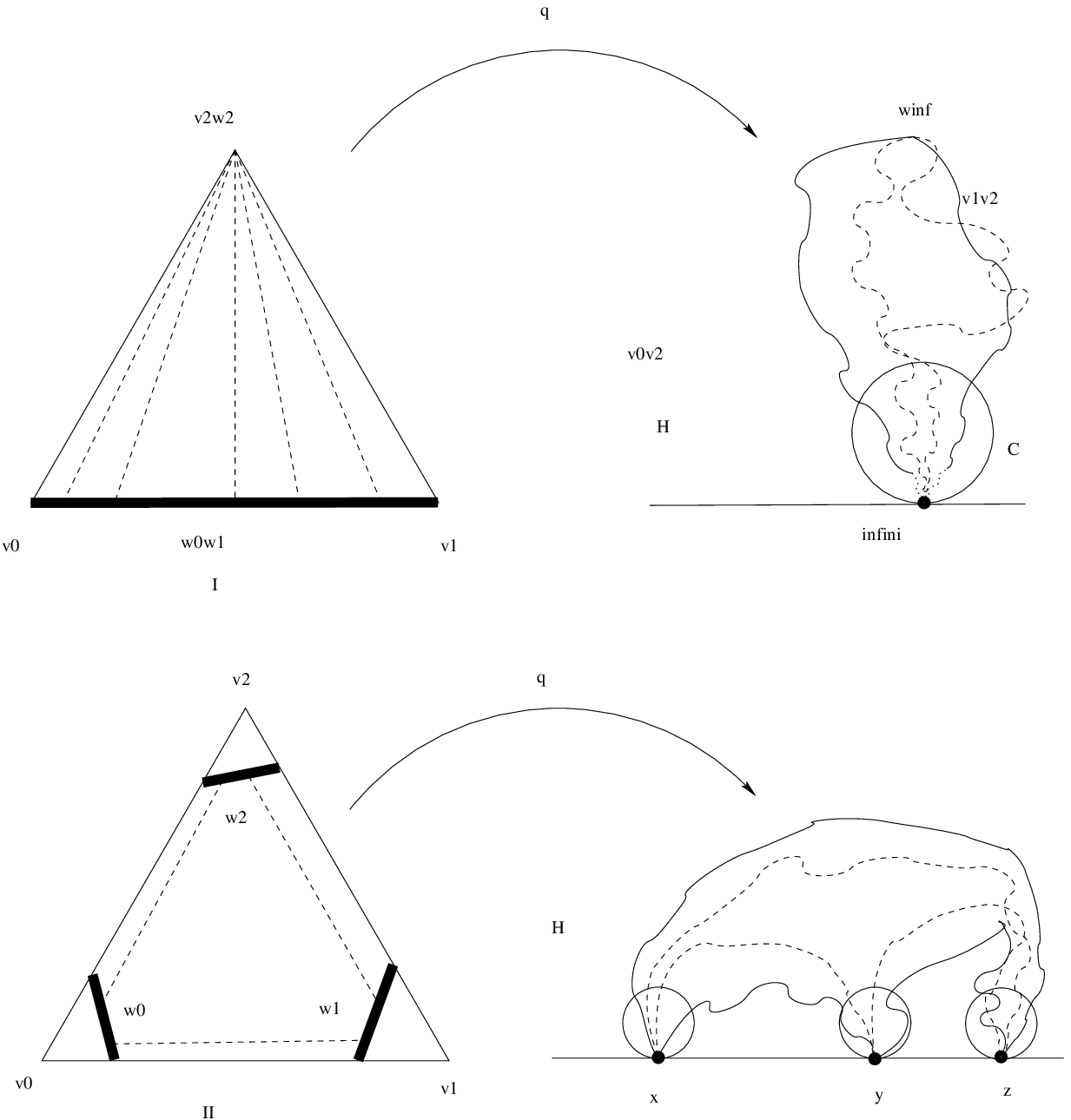}
\centerline{Figure 3}
\end{center}
\vskip 0.5 truecm
  
Let   $\{w^{\infty}_0,w^{\infty}_1,w^{\infty}_2\}$ denote three points (note necessarily pairwise distinct) of $\o{\Hi}^2$. Then we denote by $[w^{\infty}_0,w^{\infty}_1,w^{\infty}_2]$  the straight $2$-simplex of $\o{\Hi}^2$ spanned by $\{w^{\infty}_0,w^{\infty}_1,w^{\infty}_2\}$ and oriented by the order of its vertices. 
Then we associate to each component $D$ of $\sigma^{-1}({\rm int}(P))$  an algebraic area denoted by $\c{A}_{\Hi}(D)$ defined by $\pm$ the hyperbolic  area of $[w_0^{\infty},w_1^{\infty},w_2^{\infty}]$ depending on whether the orientation of $[w_0^{\infty},w_1^{\infty},w_2^{\infty}]$ matches the orientation of ${\Hi}^2$ or not. Notice that the area $\c{A}_{\Hi}(D)$ does not depend on the chosen lifting $\t{q\sigma_D}\co D\to{\Hi}^2$  of $q\sigma_D\co D\to W$ into the universal covering ${\Hi}^2$ of $W$.

 Let us define a $2$-cochain in $M$ in the following way. Let $\sigma\co\Delta^2\to M$ denote a singular $2$-simplex. If  $\rho\sigma( \Delta^2)\subset M\setminus{\rm int(P)}$ then we set 
$$\left\l\Omega_P,\sigma\right\r=0$$
and if $\rho\sigma( \Delta^2)\cap{\rm int(P)}\not=\emptyset$ then we  set 
$$\l\Omega_P,\sigma\r=\sum_{D\in\rho(\sigma)^{-1}({\rm int}(P))}\c{A}_{\Hi}(D)$$

\subsection{Proof of Proposition \ref{bounded}}\label{proofbounded}
We first check that $\Omega_P$ is a cocycle. Let $\sigma\co\Delta^3\to M$ be a $3$-simplex. Then we have by definition of $\Omega_P$
$$\l\delta\Omega_P,\sigma\r=\l\delta\Omega_P,\rho(\sigma)\r=\sum_{i=0}^3(-1)^i\l\Omega_P,\rho(\sigma)|\Delta^2_i\r$$ where $\Delta^2_i$ denotes the $i$-th face $[v_0,...,\hat{v}_i,...,v_3]$ of $\Delta^3$. To check the cocyle condition one can always assume that  $\rho\sigma(\Delta^3)\cap{\rm int}(P)\not=\emptyset$. Thus by the definition of $\Omega_P$ we have 
$$\l\delta\Omega_P,\sigma\r=\sum_{i=0}^3(-1)^i\sum_{D\in(\rho\sigma|\Delta^2_i)^{-1}({\rm int}(P))}\c{A}_{\Hi}(D)$$
Notice that for each component $D$ of $(\rho\sigma|\Delta^2_i)^{-1}({\rm int}(P))$ there exists a component $\nabla$ of  $\rho(\sigma)^{-1}({\rm int}(P))$ whose $D$ is a $2$-face. Thus to check the cocyle condition it is sufficient to prove that for any component $\nabla$ of  $\rho(\sigma)^{-1}({\rm int}(P))$ then $$\sum_{i=0}^3(-1)^i\c{A}_{\Hi}(\b_i\nabla)=0$$ where  $\b_i\nabla=\nabla\cap\Delta^2_i$ for $i=0,...,3$, with the obvious convention $\c{A}_{\Hi}(\b_i\nabla)=0$ if $\sigma(\b_i\nabla)\subset\b P$.

1. Assume that $\o{\nabla}\cap(\rho\sigma)^{-1}(\b P)$ contains a $2$-face $\Delta^2_i$ of $\Delta^3$. Then in this case each $\b_j\nabla$, for $j\not=i$, has only $2$ vertices and thus $\c{A}_{\Hi}(\b_i\nabla)=0$ for $i=0,...,3$.   

2. Assume that $\o{\nabla}\cap(\rho\sigma)^{-1}(\b P)$ contains no face  of $\Delta^3$ and two edges, say $e_1$ and $e_2$ of $\Delta^3$. Since by Remark \ref{core} we have $\o{e}_1\cap\o{e}_2=\emptyset$ then as in Case 1 each $\b_i\nabla$ has only $2$ vertices.  

3.1 Assume that $(\rho\sigma)^{-1}(\b P)\cap\o{\nabla}$ contain no face, one edge, say $[v_i,v_j]$ with $i<j$ of $\Delta^3$ and a normal rectangle $R$. Then again in this case each $\b_i\nabla$ has only $2$ vertices. 

3.2. Assume that $(\rho\sigma)^{-1}(\b P)\cap\o{\nabla}$ contain no face, one edge $[v_i,v_j]$ of $\Delta^3$ and at most two ideal vertices corresponding to $v_k$ and $v_l$.  

If $i<j<k<l$ then $\c{A}_{\Hi}(\b_{\nu}\nabla)=0$ when $\nu=k,l$ since in this case $\b_{\nu}\nabla$ has only $2$ vertices. Denote by $\{w_i,w_k,w_l\}$, resp. $\{w_j,w_k,w_l\}$, the vertices of $\b_j\nabla$, resp. of $\b_i\nabla$. Choose a base point $x$ in the open edge $(w_k,w_l)$ and a corresponding base point $\t{x}$ in ${\Hi}^2$ over $q\rho\sigma(x)$. Denote by  $\{w^{\infty}_i,w^{\infty}_k,w^{\infty}_l\}$, resp. $\{w^{\infty}_j,w^{\infty}_k,w^{\infty}_l\}$, the $\sigma$-images of  $\{w_i,w_k,w_l\}$, resp. $\{w_j,w_k,w_l\}$ corresponding to the lifting $\t{q\rho\sigma}_{|\b_j\nabla}$ of  $q\rho\sigma_{|\b_j\nabla}$, resp.  $\t{q\rho\sigma}_{|\b_i\nabla}$ of  $q\rho\sigma_{|\b_i\nabla}$, such that $\t{q\rho\sigma}_{|\b_j\nabla}(x)=\t{q\rho\sigma}_{|\b_i\nabla}(x)=\t{x}$. Notice that since $\rho\sigma([v_i,v_j])\subset\b P$ then $w_i^{\infty}=w_j^{\infty}$. Thus in this case, since $j=i+1$, we have the equality 
$$(-1)^j\c{A}_{\Hi}(\b_j\nabla)+(-1)^i\c{A}_{\Hi}(\b_i\nabla)=0$$  

If $i<k<j<l$ then $\c{A}_{\Hi}(\b_{\nu}\nabla)=0$ when $\nu=k,l$. Denote by $\{w_i,w_k,w_l\}$, resp. $\{w_k,w_j,w_l\}$, the vertices of $\b_j\nabla$, resp. of $\b_i\nabla$. Denote by  $\{w^{\infty}_i,w^{\infty}_k,w^{\infty}_l\}$, resp. $\{w^{\infty}_k,w^{\infty}_j,w^{\infty}_l\}$, the $\sigma$-images of  $\{w_i,w_k,w_l\}$, resp. $\{w_k,w_j,w_l\}$. Since $\rho\sigma([v_i,v_j])\subset\b P$ then $w_i^{\infty}=w_j^{\infty}$. Thus in this case we have the equality 
$$(-1)^j\c{A}_{\Hi}(\b_j\nabla)+(-1)^i\c{A}_{\Hi}(\b_i\nabla)=0$$  since $j=i+2$ and $[w^{\infty}_i,w^{\infty}_k,w^{\infty}_l]$ and  $[w^{\infty}_k,w^{\infty}_j,w^{\infty}_l]$ are the same geodesic simplices with opposite orientations. The other possibilities for $i,j,k,l$ follow in the same way.

Assume now that $(\rho\sigma)^{-1}(\b P)\cap\o{\nabla}$ consists only of vertices, normal triangles or normal rectangles of  $\Delta^3$.  Note that it follows from the construction that the only possibilities are ${\rm Card}(V(\nabla))=0,2,4$. 

4.1.  Assume ${\rm Card}(V(\nabla))= 4$ and denote by $w_0,w_1,w_2,w_3$ the vertices of $\nabla$ corresponding to $v_0,v_1,v_2,v_3$ and denote by $w_0^{\infty},w_1^{\infty},w_2^{\infty},w_3^{\infty}$ their $\sigma$-images. Then  
$$\sum_{i=0}^3(-1)^i[w_0^{\infty},...,\hat{w_i^{\infty}},...,w_3^{\infty}]$$ is the boundary of a geodesic $3$-simplex in ${\Hi}^2$. Using the same argument as in \cite{BG} we know that there are two distinct configurations for a geodesic $3$-simplex in ${\Hi}^2$ and in any case we have 
 $$\sum_{i=0}^3(-1)^i\c{A}_{\Hi}({\b_i\nabla})=0$$

4.2. Assume ${\rm Card}(V(\nabla))=2$. Denote by $w_i$ and $w_j$ the vertices of $\nabla$.  Then we have the two following possibilities:

4.2.1. Suppose first that $\b\nabla$ contains a normal rectangle $R$. We denote by $v_i, v_j$ the vertices of $\Delta^3$ corresponding to the vertices $w_i$ and $w_j$ of $\nabla$. First notice that $\b_i\nabla$ and $\b_j\nabla$ contain only two vertices: $w_j$, resp. $w_i$, and the vertex corresponding to $R\cap\b_i\nabla$, resp. $R\cap\b_j\nabla$. Then 
$$\c{A}_{\Hi}({\b_i\nabla})=\c{A}_{\Hi}({\b_j\nabla})=0$$
On the other hand $\b_k\nabla$ and $\b_l\nabla$ have three vertices: the vertices $w_i$ and $w_j$ and one ideal vertex (denoted by $w_l$, resp. $w_k$) "corresponding" to $R$. Since $w_l$ and $w_k$ correspond to edges of the same rectangle then $w_k^{\infty}=w_l^{\infty}$ (for some  lifting with the same base point). Then we get again as in 3.2

$$(-1)^k\c{A}_{\Hi}(\b_k\nabla)+(-1)^l\c{A}_{\Hi}(\b_l\nabla)=0$$

4.2.2 Suppose that $\b\nabla$ contains no normal rectangle $R$ then  

 either $V(\nabla)$ consists of one normal triangle and one vertex of $\Delta^3$. Then  $\c{A}_{\Hi}(\b_i\nabla)=0$ for any $i$, since $\b_i\nabla$ contains only two vertices. 

 either $V(\nabla)$ consists of two parallel normal triangles. Then again $\c{A}_{\Hi}(\b_i\nabla)=0$ for any $i$, since $\b_i\nabla$ contains only two vertices. 
 
4.3. Assume ${\rm Card}(V(\nabla))=0$. Then necessarily, $\b\nabla$ contains two parallel rectangles and $\c{A}_{\Hi}(\b_i\nabla)=0$ for any $i$ since $\b_i\nabla$ contains only two vertices. This completes the  proof that $\Omega_P$ is a cocyle. 

Next point (i) of the proposition follows directly from the construction of $\Omega_P$. On the other hand given $\alpha\in\hat{H}_2(P,\b P;{\R})$ and $z_{\alpha}$ a relative $2$-cycle representing of $\alpha$, we get the following equalities
$$\left\l[\Omega_P],\hat{\Theta}_P(\alpha)\right\r=\left\l\Omega_P,z_{\alpha}+u_{\alpha}\right\r=\left\l i^{\ast}(\Omega_P),z_{\alpha}\right\r=\left\l [i^{\ast}(\Omega_P)],[z_{\alpha}]\right\r$$
where $u_{\alpha}$ is a $l_1$-chain in $\b P$ such that $\b u_{\alpha}=-\b z_{\alpha}$. In particular, if $z_{\alpha}$ represents the fundamental class of $F$ then we get $$\left\l[\Omega_P],\hat{\Theta}_P(\alpha)\right\r={\rm Area}(F)$$ Indeed we can choose as a representant $z_{\alpha}$ the formal sum of a triangulation of $F$ and then apply $\Omega_P$. On the other hand, if $\c{F}$ denotes any connected horizontal surface in $P$ then since $\c{F}$ is a finite covering of $F$ we get  $\left\l[\Omega_P],\hat{\Theta}_P([\c{F}])\right\r={\rm Area}(\c{F})$. This proves point (ii). 

It remains to check point (iii). The fact that $k_i^{\ast}(\Omega_{P_j})=0$ for any $i\not=j$ follows from the construction of $\Omega_{P_j}$.  
We first check that $\|\sum_{i\in I}\Omega_{P_i}\|_{\infty}\leq\pi$. For each $i\in I$ we can identify $P_i$ with the product $F_i\times{\S}^1$. Let $\sigma\co\Delta^2\to M$ be a singular $2$-simplex. If there exists an edge $e$ of $\Delta^2$ such that $\rho\sigma(e)\subset\c{T}$ then $\l\sum_{i\in I}\Omega_{P_i},\sigma\r=0$. If for any edge $e$ of $\Delta^2$ we have  $\rho\sigma(e)\not\subset\c{T}$ then we know from Remark \ref{core} that there exists a unique component ${\rm Core}(\sigma)$ of $(\rho\sigma)^{-1}(M\setminus\c{T})$ which meets the three edges of $\Delta^2$.  Denote by $P_{\nu}$ the component of $M\setminus\c{T}$ such that $\rho\sigma({\rm Core}(\sigma))\subset{\rm int}(P_{\nu})$. If $D$ is a component of $(\rho\sigma)^{-1}({\rm int}(P_i))$ for some $i\in I$ distinct from ${\rm Core}(\sigma)$ then we get the equality

$$\c{A}_{\Hi}(D)=0$$
 Indeed in this case  $D$ has $2$ (ideal) vertices.  On the other hand, if $P_{\nu}=P_{i_0}$ for some $i_0\in I$ then   we know that  $$\left|\c{A}_{\Hi}(D)\right|\leq\pi$$ by the definition of $\c{A}_{\Hi}(D)$ since the areas of geodesic triangles in ${\Hi}^2$ are bounded by $\pi$. This proves that $\|\sum_{i\in I}\Omega_{P_i}\|_{\infty}\leq\pi$. 

Denote by $\Omega$ the sum $\sum_{i\in I}\Omega_{P_i}$ and let $i\in I$ be a fixed index. Since $\Theta_{P_i}$ is a contraction, then using point (ii)  we get the following equalities
$$\left|\left\l[\Omega],\hat{\Theta}_{P_i}([F_i])\right\r\right|={\rm Area}(F_i)\leq\left\|[\Omega]\right\|_{\infty}\left\|\hat{\Theta}_{P_i}([F_i])\right\|_1\leq\left\|[\Omega]\right\|_{\infty}\left\|[F_i]\right\|_1$$ Since $\|[F_i]\|_1\leq\|F_i\|$
this completes the proof of the proposition since ${\rm Area}(F_i)=\pi\|F_i\|$ by \cite[Section 0.4]{G} and \cite{Th}.
\section{A Filling isometry for relative classes}\label{filling}
This section is devoted to the proof of Theorem \ref{1} and of Proposition \ref{2}. We begin with a special case.
\begin{lemma}\label{1facile}
  Theorem \ref{1} is true when each Seifert piece of $(M,\c{T})$ is homeomorphic to a product $F\times{\S}^1$.
\end{lemma}   
\begin{proof}
We first check point (i). Let $P$ be a Seifert piece of $M\setminus\c{T}$. Then $P$ is homeomorphic to a product $F\times{\S}^1$ where $F$ is a hyperbolic orientable surface. By Proposition \ref{bounded} we deduce, using the Holder inequality that $\|\hat{\Theta}_P([F])\|_1\geq\|F\|\geq\|[F]\|_1$ and thus finally we get from above $\|\hat{\Theta}_P([F])\|_1=\|[F]\|_1=\|F\|$.  Let $\alpha\in\hat{H}_2(P,\b P;{\R})$. Then by Lemma \ref{generator}, there exists $\xi\in{\R}$ such that $\alpha=\xi\alpha_P$. Then using Proposition \ref{bounded} we get
$$\left|\left\l[\Omega_P],\hat{\Theta}_P(\alpha)\right\r\right|=|\xi|{\rm Area}(F)\leq\pi\|\hat{\Theta}_P(\alpha)\|_1$$
Hence $\|\hat{\Theta}_{P}(\alpha)\|_1\geq|\xi|\|F\|\geq\|\alpha\|_1$ and thus $\hat{\Theta}_P$ is an isometry. If $\alpha$ is the class of horizontal surface $\c{F}$ then $F$ is finitely covered by $\c{F}$  and thus  $\xi$  is necessarily and integer satisfying  $\|\c{F}\|=|\xi|\|F\|$ (actually  $\xi$ is the degree of the map induced by the Seifert projection $\c{F}\to F$). This proves point (i).

We check point (ii). Let $P_1,...,P_l$ denote the Seifert pieces of $M\setminus\c{T}$. 
 For each $P_i$ denote by $F_i$ a minimal surface. Let $(\alpha_1,...\alpha_l)$ be a $l$-uple of $\hat{H}_2(P_1,\b P_1)\times...\times\hat{H}_2(P_l,\b P_l)$ and  denote by $\xi_1,...,\xi_l$ the real numbers such that for each $i=1,...,l$ we have $\alpha_i=\xi_i\alpha_{P_i}$, where $\alpha_{P_i}$ denotes the class of $F_i$ in $\hat{H}^+_2(P_i,\b P_i)$. By Proposition \ref{bounded} there exists for each $i\in\{1,...,l\}$ a bounded $2$-cocycle $\Omega_i$ in $M$ such that $$\left\l[\Omega_i],\hat{\Theta}_{P_j}\alpha_j\right\r=\delta_{ij}|\xi_i|{\rm Area}(F_i)$$
 for any $i,j$ in $\{1,...,l\}$. Thus we get   
$$\left\l\sum_i[\Omega_{i}],\sum_j\hat{\Theta}_{P_j}(\alpha_j)\right\r=\sum_i|\xi_i|{\rm Area}(F_i)\leq\pi\left\|\sum_j\hat{\Theta}_{P_j}(\alpha_{j})\right\|_1$$
thus we get 
$$\left\|\sum_i\hat{\Theta}_{P_i}(\alpha_i)\right\|_1\geq\sum_i|\xi_i|\|F_i\|\geq\sum_i\left\|\xi_i[F_i]\right\|_1=\sum_i\left\|\hat{\Theta}_{P_i}(\alpha_{i})\right\|_1$$
Since the opposite inequality is also true this proves the lemma.
\end{proof}
We now establish the following result which is a special case of Proposition \ref{2}.
\begin{lemma}\label{galoisien}
Let $(M,\c{T})$ be a closed aspherical framed $3$-manifold endowed with a torus splitting and let $p\co\t{M}\to M$ denote a finite regular covering such that    each Seifert component of $\t{M}\setminus\t{\c{T}}$  is homeomorphic to a product, where $\t{\c{T}}=p^{-1}(\c{T})$. Then the covering $p$ induces an isometry $p_{\sharp}|\hat{H}_{2}^{l_1,+}(\t{M}^{\t{\c{T}}};{\R})\co\hat{H}_2^{l_1,+}(\t{M}^{\t{\c{T}}};{\R})\to \hat{H}_2^{l_1}(M;{\R})$ where $\t{M}$ is endowed with the framing induces by that of $M$.
\end{lemma}
\begin{proof}  Denote by $\Gamma$ the automorphism group of $p\co\t{M}\to M$ and by $\{\t{P}_1,...,\t{P}_l\}$ the Seifert pieces of $(\t{M},\t{\c{T}})$. For each $i=1,...,l$ we know that $\t{P}_i$ is homeomorphic to a product $\t{F}_i\times{\S}^1$.  Let $\t{\alpha}$ be an element of $\hat{H}_2^{l_1,+}(\t{M}^{\t{\c{T}}};{\R})$ and denote by ${\rm Av}(\t{\alpha})$ the class obtained by \emph{averaging}  $\t{\alpha}$ defined by 
$${\rm Av}(\t{\alpha})=\sum_{g\in\Gamma}g_{\sharp}(\t{\alpha})\in\Gamma\hat{H}_2^{l_1}(\t{M};{\R})$$
By definition we have $\t{\alpha}=\sum_{i\in I}\xi_i\hat{\Theta}_{\t{P}_i}(\alpha_{\t{P}_i})$ where $I$ is a subset of $\{1,...,l\}$ and $\xi_i\in{\R}_+$ and thus since each $g\in\Gamma$ preserves the torus splitting by the naturality property of Lemma \ref{rempli} we get 
$${\rm Av}(\t{\alpha})=\sum_{g\in\Gamma,i\in I}\xi_i\hat{\Theta}_{g(\t{P}_i)}\left(g_{\sharp}\left(\alpha_{\t{P}_i}\right)\right)$$
Moreover notice that $$\left\|\hat{\Theta}_{g(\t{P}_i)}g_{\sharp}\alpha_{\t{P}_i}\right\|_1\leq\left\|\t{F}_i\right\|$$


For each $i=1,...,l$ denote by $\Omega_{\t{P}_i}$  the bounded $2$-cocycle of $\t{M}$ constructed in Proposition \ref{bounded} and denote by $\Omega$ the sum $\sum_{i=1}^l\Omega_{\t{P}_i}$.  Since  each $g\in\Gamma$ acts one $\t{M}$ as an orientation preserving diffeomorphism which preserves the orientation of the fibers (by the framing hypothesis)  and the torus decomposition $\t{\c{T}}$ of $\t{M}$ then we get
$$\left\l[\Omega],{\rm Av}(\t{\alpha})\right\r={\rm Card}(\Gamma)\sum_{i\in I}\xi_i{\rm Area}(\t{F}_i)\leq\pi\|{\rm Av}(\t{\alpha})\|_1$$
This proves that $$\left\|{\rm Av}(\t{\alpha})\right\|_1={\rm Card}(\Gamma)\sum_{i\in I}\xi_i\left\|\t{F}_i\right\|$$ 
Recall that
since ${\rm Av}(\t{\alpha})\in\Gamma H_2^{l_1}(\t{M};{\R})$ then by Lemma \ref{galois} $\|p_{\sharp}({\rm Av}(\t{\alpha}))\|_1=\|{\rm Av}(\t{\alpha})\|_1$.  On the other hand we have
$$\left\|p_{\sharp}({\rm Av}(\t{\alpha}))\right\|_1\leq\sum_{g\in\Gamma}\left\|p_{\sharp}g_{\sharp}(\t{\alpha})\right\|_1\leq\sum_{g\in\Gamma}\left\|g_{\sharp}(\t{\alpha})\right\|_1\leq{\rm Card}(\Gamma)\sum_{i\in I}\xi_i\left\|\t{F}_i\right\|$$
We deduce that $\sum_{g\in\Gamma}\left\|p_{\sharp}g_{\sharp}(\t{\alpha})\right\|_1=\sum_{g\in\Gamma}\left\|g_{\sharp}(\t{\alpha})\right\|_1$. On the other hand, since we know that $\left\|p_{\sharp}g_{\sharp}(\t{\alpha})\right\|_1\leq\left\|g_{\sharp}(\t{\alpha})\right\|_1$ for any $g\in\Gamma$ then we get in particular $\left\|p_{\sharp}(\t{\alpha})\right\|_1=\left\|\t{\alpha}\right\|_1$.

\end{proof}

  \begin{proof}[Proof of Theorem \ref{1}] Let $p\co\t{M}\to M$ be a finite regular covering of $M$ endowed with a torus splitting defined by $\t{\c{T}}=p^{-1}(\c{T})$ and such that   each component of  $\t{M}\setminus\t{\c{T}}$  is a product (such a covering exists by \cite[Proposition 4.4]{LW}). For each Seifert piece $P$ of $M\setminus\c{T}$ denote by $\t{P}$ a component over $P$ in $\t{M}$. 
Let $\alpha\in\hat{H}_2(P,\b P;{\R})$. Then there exists $\t{\alpha}\in\hat{H}_2(\t{P},\b\t{P};{\R})$ such that $p_{\sharp}(\t{\alpha})=\alpha$. Hence we have, using Lemmas \ref{1facile} and  \ref{galoisien}, $\|\hat{\Theta}_P(\alpha)\|_1=\|\hat{\Theta}_P(p_{\sharp}(\t{\alpha}))\|_1=\|p_{\sharp}\hat{\Theta}_{\t{P}}(\t{\alpha})\|_1=\|\hat{\Theta}_{\t{P}}(\t{\alpha})\|_1=\|\t{\alpha}\|_1\geq\|\alpha\|_1$. Indeed replacing $\t{\alpha}$ by $-\t{\alpha}$ we may assume that $\hat{\Theta}_{\t{P}}(\t{\alpha})\in\hat{H}_2^{l_1,+}(\t{M}^{\t{\c{T}}})$.  Since $\hat{\Theta}_P$ is a contraction by Lemma \ref{rempli} this proves the isometry. 
To complete the proof of point (i) it is sufficient to check that if $\alpha$ is the class of a connected horizontal surface $\c{F}$ in $P$ then $\|\hat{\Theta}_P(\alpha)\|_1=\|\c{F}\|$. Denote by $\t{\c{F}}$ a component over $\c{F}$ in $\t{P}\subset\t{M}$ and denote by $n$ the degree of the induced covering map $\t{\c{F}}\to\c{F}$.  
Let $\t{\alpha}\in\hat{H}_2(\t{P},\b\t{P})$ be the integral class corresponding to $\t{\c{F}}$.  Then we know from Lemma \ref{1facile} that $\|\hat{\Theta}_{\t{P}}(\t{\alpha})\|_1=\|\t{\c{F}}\|=\|\t{\alpha}\|_1$. On the other hand we know  that 
$p_{\sharp}(\hat{\Theta}_{\t{P}}(\t{\alpha}))=n\hat{\Theta}_P(\alpha)$. Then we get using Lemma \ref{galoisien}
$$\left\|\t{\c{F}}\right\|=\left\|\t{\alpha}\right\|_1=\left\|\hat{\Theta}_{\t{P}}(\t{\alpha})\right\|_1=\left\|p_{\sharp}(\hat{\Theta}_{\t{P}}(\t{\alpha}))\right\|_1=|n|\left\|\hat{\Theta}_P(\alpha)\right\|_1\leq|n|\|\alpha\|_1$$ Since $\|\alpha\|_1\leq\|\c{F}\|$
and $\|\t{\c{F}}\|=n\|\c{F}\|$ then we get  $\|\hat{\Theta}_P(\alpha)\|_1=\|\alpha\|_1=\|\c{F}\|$. This completes the proof of point (i).

We check the additivity property of the $l_1$-norm.
 For each $P_i$ denote by $F_i$ a minimal surface  and for each $i=1,...,l$ choose an element $\alpha_i\in\hat{H}^+_2(P_i,\b P_i;{\R})$ and we denote by $\alpha$ the $l_1$-class given by $\sum_{i=1}^l\hat{\Theta}_{P_i}(\alpha_i)$. Denote by $\xi_1,...,\xi_l$ the non-negative real numbers such that for each $i=1,...,l$ we have $\alpha_i=\xi_i\alpha_{P_i}$, where $\alpha_{P_i}$ denotes the class of $F_i$ in $\hat{H}_2(P_i,\b P_i)$. For each $i$ there exists an element $\t{\alpha}_i$ in $\hat{H}^+_2(\t{P}_i,\b \t{P}_i;{\R})$ such that $p_{\sharp}(\t{\alpha}_i)=\alpha_i$. Denote by $\t{\alpha}$ the element $\sum_{i=1}^l\hat{\Theta}_{\t{P}_i}(\t{\alpha}_i)\in\hat{H}_2^{l_1,+}(\t{M}^{\t{\c{T}}};{\R})$ such that $p_{\sharp}(\t{\alpha})=\alpha$. Since we may also assume that $\hat{\Theta}_{\t{P}_i}(\t{\alpha}_i)\in\hat{H}_2^{l_1,+}(\t{M}^{\t{\c{T}}})$ for each $i=1,...,l$ it follows from Lemma \ref{galoisien} combined with Lemma \ref{1facile} that 
$$\|\alpha\|_1=\|\t{\alpha}\|_1=\sum_i\|\hat{\Theta}_{\t{P}_i}(\t{\alpha}_i)\|_1=\sum_i\|p_{\sharp}\hat{\Theta}_{\t{P}_i}(\t{\alpha}_i)\|_1=\sum_i\|\hat{\Theta}_{P_i}(\alpha_i)\|_1$$  
This completes the proof of Theorem \ref{1}.
\end{proof}
\begin{proof}[Proof of Proposition \ref{2}]
 Let $p\co\t{M}\to M$ denote a finite covering and let $q\co\hat{M}\to\t{M}$ be a finite covering endowed with a torus splitting defined by $\hat{\c{T}}=q^{-1}(\t{\c{T}})$ such that each Seifert piece of $(\hat{M},\hat{\c{T}})$ is homeomorphic to a product and such that $r=p\circ q$ is regular. Let $\{\t{P}_1,...,\t{P}_l\}$ denote a collection of Seifert pieces of $(\t{M},\t{\c{T}})$ and for each $i$ denote by $\t{\alpha}_i$ an element of  $\hat{H}^+_2(\t{P}_i,\b \t{P}_i;{\R})$.    Let $\t{\alpha}$ denote the $l_1$-class equal to $\sum_{i=1}^l\hat{\Theta}_{\t{P}_i}(\t{\alpha}_i)$. Then using the above construction we know that there exists $\hat{\alpha}\in \hat{H}^{l_1,+}_2(\hat{N}^{\hat{\c{T}}};{\R})$ such that $q_{\sharp}(\hat{\alpha})=\t{\alpha}$ and $\|\hat{\alpha}\|_1=\|r_{\sharp}(\hat{\alpha})\|_1$ by Lemma \ref{galoisien}. Thus we get $\|p_{\sharp}\t{\alpha}\|_1=\|p_{\sharp}q_{\sharp}\hat{\alpha}\|_1=\|r_{\sharp}\hat{\alpha}\|_1=\|\hat{\alpha}\|_1\geq\|\t{\alpha}\|_1$. Since any continuous map induces a contraction between $l_1$-homology groups we deduce that $\|p_{\sharp}\t{\alpha}\|_1=\|\t{\alpha}\|_1$ which completes the proof of the proposition.
\end{proof} 

\section{Characterizations of finite covering maps}

We begin this section with some recall on the Jaco Shalen Johannson torus decomposition of $3$-manifolds which will be used throughout  the proof of Theorem \ref{localiso} and Corollary \ref{rigidity}.
Given a closed irreducible orientable $3$-manifold $N$ we denote by ${\c{T}}_N$ the Jaco-Shalen-Johannson family of canonical tori of $N$ and by ${\c{H}}(N)$ (resp. ${\c{S}}(N)$) the disjoint union of the hyperbolic (resp. Seifert) components of $N\setminus {\c{T}}_N\times(-1,1)$ so that $N\setminus{\c{T}}_N\times(-1,1)={\c{H}}(N)\cup{\c{S}}(N)$, where ${\c{T}}_N\times[-1,1]$ is identified with a regular neighborhood of ${\c{T}}_N$ in such a way that ${\c{T}}_N\simeq{\c{T}}_N\times\{0\}$  (see \cite{JS}, \cite{J} and \cite{Th1}). On the other hand, we denote by $\Sigma(N)=(\Sigma(N),\emptyset)$ the \emph{characteristic Seifert pair of $N$} in the sense of \cite{JS} and \cite{J}. We start by recalling a main consequence of the Characteristic Pair Theorem of W. Jaco and P. Shalen (see \cite[Chapter V]{JS}) which allows to control a nondegenerate map from a Seifert fibered space into an  irreducible $3$-manifold.  We first give the definition of degenerate maps in the sense of W. Jaco and P. Shalen.
\begin{definition}\label{deg}
Let $(S,F)$ be a connected Seifert pair, and let $(N,T)$ be a connected 3-manifold pair. A map $f\co(S,F)\to(N,T)$ is said to be \emph{degenerate} if either

(0) the map $f$ is inessential as a map of pairs (i.e. $f$ is homotopic, as a map of pairs, to a map $g$ such that $g(S)\subset T$), or

(1) the group ${\rm Im}(f_{\ast}\co\pi_1S\to\pi_1N)=\{1\}$, or

(2) the group ${\rm Im}(f_{\ast}\co\pi_1S\to\pi_1N)$ is cyclic and $F=\emptyset$, or

(3) the map $f|\gamma$ is homotopic in $N$ to a constant map for some fiber $\gamma$ of $(S,F)$.    
\end{definition}
Then the Characteristic Pair Theorem of Jaco and Shalen implies the following result.
\begin{theorem}\label{ndeg}[Jaco, Shalen]
If $f$ is a nondegenerate map of a Seifert pair $(S,\emptyset)$ into a  closed irreducible orientable $3$-manifold $(N,\emptyset)$, then there exists a map $f_1$ of $S$ into $N$, homotopic to $f$, such that $f_1(S)\subset{\rm int}(\Sigma(N))$.
\end{theorem}

In order  to prove Theorem \ref{localiso} we first check the following technical result. 

\begin{proposition}\label{reduction} Let $M$ be a closed aspherical oriented $3$-manifold.  
Any $\pi_1$-surjective nonzero degree map $f\co M\to N$  into a  closed irreducible orientable $3$-manifold satisfying the following conditions

(i) Each Seifert component of $M\setminus\c{T}_M$, resp. of $N\setminus\c{T}_N$, is homeomorphic to a product, resp. to a ${\S}^1$-bundle over an orientable surface,  each Seifert component of  $M\setminus\c{T}_M$  has at least two boundary components (if $\c{T}_M\not=\emptyset$) and each component of $\c{T}_M$ is shared by two distinct components of $M\setminus\c{T}_M$,

(ii) $\|f_{\sharp}[M]\|_1=\|[M]\|_1$, where $[M]\in H_3(M;{\R})$ is the fundamental class

(iii)   $\|f_{\sharp}\hat{\Theta}_P({\alpha})\|_1=\|\hat{\Theta}_P\alpha\|_1$ for each $\alpha\in\hat{H}_2(P,\b P)$ when $P$ runs over the Seifert components of $M\setminus\c{T}_M$

 is homotopic to a homeomorphism.
 
\end{proposition}

\subsection{Proof of Proposition \ref{reduction}}

Throughout this section we always assume that the map $f\co M\to N$ and the manifolds $M, N$ satisfy the hypothesis of Proposition \ref{reduction}. Notice that we may assume in addition  that $M$ is not a virtual torus bundle by \cite{W}. Thus since each Seifert piece $P$ of $M$ is homeomorphic to a product we may assume that $P$ is a ${\Hi}^2\times{\R}$-manifold.  Hence this implies, using hypothesis (ii) and (iii), that either $\|N\|\not=0$ or $\hat{H}_2^{l_1}(N;{\R})\not=\{0\}$. Hence either $N$ has a non-empty JSJ-splitting or $N$ admits a geometry ${\Hi}^3, {\Hi}^2\times{\R}$ or $\t{\rm SL}(2,{\R})$.   The proof of Proposition \ref{reduction} will come from the following sequence of claims.

\begin{claim}\label{surface}
Let  $P$ be an incompressible ${\Hi}^2\times{\R}$-submanifold of $M$ and let $c$ be a simple closed curve in some component $T$ of $\b P$. If there exists a horizontal surface in $P$ whose $c$ is a boundary component then $f_{\ast}([c])\not=1$ in $\pi_1N$. Moreover, $f_{\ast}(\pi_1P)$ is a non-abelian group.
\end{claim}
\begin{proof}
Denote by $F$ a horizontal surface in $P$ whose $c$ is a boundary component. Since $P$ is a ${\Hi}^2\times{\R}$-manifold then $F$ is necessarily a hyperbolic surface. In particular $F$ has a positive simplicial volume. Suppose that $f_{\ast}([c])=1$. Denote by $T\times[-1,1]$ a regular neighborhood of $T$ such that $T=T\times\{0\}$ and parametrize $T={\S}^1\times{\S}^1$ such that $c={\S}^1\times\{\ast\}$. As in \cite{Ro}, consider the relation $\sim$ on $M$ defined by $z\sim z'$ iff $z=z'$ or $z=(x,y,t)\in T\times I$, $z'=(x',y',t')\in T\times I$ and $y=y'$, $t=t'$. Denote by $X=M/\sim$ the quotient space and by $\pi\co M\to X$ the quotient map. Then the map $f$ factors throught $X$. Denote by $g\co X\to N$ the map such that $f=g\circ\pi$. Denote by $\hat{P}$ the image of $P$ under $\pi$. Topologically $\hat{P}$ is obtained from $P$ after Dehn filling along $T$, identifying the meridian of a solid torus $V$ to $c$ so that the Seifert fibration of $P$ extends to a Seifert fibration of $\hat{P}$ and the image $\hat{F}$ of $F$ is a  surface in $\hat{P}$ obtained from $F$ after gluing a $2$-disk along each component of $\b F$ parallel to  $c$. Denote by $\c{C}$ the union of the components of $\b F$ parallel to $c$. Note that it follows from our construction that 
$$\pi_{\sharp}(\hat{\Theta}_P([F]))=\hat{\Theta}_{\hat{P}}([\hat{F}])\in H_2^{l_1}(\hat{P};{\R})$$
Indeed denote by $z_F$ a relative cycle in $F$ representing the fundamental class of $F$ and denote by $u_F$ a $l_1$-chain in $\b F$ such that $\b u_F=-\b z_F$. Then have the following equalities: $\pi_{\sharp}(\hat{\Theta}_P([F]))=\pi_{\sharp}([z_F+u_F])=[\pi_{\sharp}(z_F)+\pi_{\sharp}(u_F)]\in\hat{H}_2^{l_1}(\hat{P})$. Next consider the induces homomorphism 
 $$\pi_{\sharp}\co H_2(F,\b F)\to H_2(F/\c{C},\b F/\c{C})=H_2(\hat{F},\b\hat{F}\coprod X_0)$$ where $X_0=\pi(\c{C})=\c{C}/\c{C}$. It follows easily from the Excision Theorem that $\pi_{\sharp}$ is actually an isomorphism so that $\pi_{\sharp}(z_F)$ represents a generator of $H_2(\hat{F},\b\hat{F}\coprod X_0)$. On the other hand, there exists a $2$-chain $u_0$ in $X_0$ such that $\pi_{\sharp}(z_F)+u_0$ is a relative cycle in $(\hat{F},\b\hat{F})$. Note that $[\pi_{\sharp}(z_F)+u_0]=[\pi_{\sharp}(z_F)]$ in $H_2(\hat{F},\b\hat{F}\coprod X_0)$. Since the inclusion $(\hat{F},\b\hat{F})\hookrightarrow(\hat{F},\b\hat{F}\coprod X_0)$ induces an isomorphism $H_2(\hat{F},\b\hat{F})\to H_2(\hat{F},\b\hat{F}\coprod X_0)$ then  $\pi_{\sharp}(z_F)+u_0$ represents a generator of $H_2(\hat{F},\b\hat{F})$. Denote by $\hat{u}_F$ a $l_1$-chain in $\b\hat{F}$ such that $\b\hat{u}_F=-\b(\pi_{\sharp}(z_F)+u_0)$. Then since each component of $\pi(\b P)$ has an amenable fundamental group then $[\pi_{\sharp}(z_F)+\pi_{\sharp}(u_F)]=[\pi_{\sharp}(z_F)+u_0+\hat{u}_F]$ in $H_2^{l_1}(\hat{P};{\R})$. 
 
 Since $\pi_{\sharp}(z_F)+u_0$ represents a generator of $H_2(\hat{F},\b\hat{F})$ then it follows that $\pi_{\sharp}(\hat{\Theta}_P([F]))=\hat{\Theta}_{\hat{P}}([\hat{F}])\in H_2^{l_1}(\hat{P};{\R})$.
 Finally we deduce the following equalities:
$$\|\hat{\Theta}_P([F])\|_1\geq\|\pi_{\sharp}\hat{\Theta}_P([F])\|_1=\|\hat{\Theta}_{\hat{P}}([\hat{F}])\|_1\geq\|f_{\sharp}\hat{\Theta}_P([F])\|_1=\|\hat{\Theta}_P([F])\|_1$$
Thus we get the following equalities:

$$\|\hat{F}\|\geq\|\hat{\Theta}_{\hat{P}}([\hat{F}])\|_1=\|\hat{\Theta}_P([F])\|_1=\|F\|$$
A contradiction. This proves the first statement of the lemma. 

It remains to check that  $f_{\ast}(\pi_1P)$ is a non-abelian group.  Assume that $f_{\ast}(\pi_1P)$ is an abelian subgroup of $\pi_1N$. Then necessarily $f_{\ast}(\pi_1P)$ is isomorphic to a free abelian group of rank $\leq 3$. Denote by $X$ a $K(f_{\ast}(\pi_1P),1)$-space ($X$ is homeomorphic to either ${\bf D}^3,{\S}^1\times{\bf D}^2, {\S}^1\times{\S}^1\times I$ or ${\S}^1\times{\S}^1\times{\S}^1$). Then the map $f|P\co P\to N$ factors through $X$ so that we have the following commutative diagram
$$\xymatrix{
\hat{H}_2(P,\b P) \ar[d]^{\simeq} \ar[r]^{\hat{\Theta}_P} & \hat{H}_2^{l_1}(P) \ar[r]^{i_{\sharp}} \ar[rd] & \hat{H}_2^{l_1}(M) \ar[r]^{f_{\sharp}} & \hat{H}_2^{l_1}(N) \\
{\R} & &  \hat{H}_2^{l_1}(X) \ar[ur]}$$ 
 
 Since  $\hat{H}_2^{l_1}(X)$ is trivial (because $X$ has an amenable fundamental group) then we get a contradiction with hypothesis (iii) of  Proposition \ref{reduction}.
\end{proof}

\begin{claim}\label{primitive}
Let $P=F\times{S}^1$ be a $3$-manifold such that $F$ is an orientable hyperbolic surface with non connected boundary. Then for any simple closed curve $c$ of $\b P$ that is not homotopic to the fiber of $P$, there exists a horizontal surface $(H,\b H)$ in $(P,\b P)$ such that $c$ is parallel to a component of $\b H$.
\end{claim}
\begin{proof} Denote by $T_1$ the component of $\b P$ which contains $c$ and by 
 $T_2,...,T_r$ the other components of $\b P$ with $r\geq 2$. For each $i=1,...,r$ fix a basis $\l s_i,h\r$, where $s_i$ is a section of $T_i$ with respect to the ${\S}^1$-fibration of $P$ such that $s_1+...+s_r$ is nul-homologous in $P$ and where $h$ denotes  the  fiber of $P$. Denote by $(\alpha,\beta)$ the coprime integers with $\alpha\not=0$ such that $c=\alpha[s_1]+\beta[h]$.  Then 
 
 $$[c]+ \alpha[s_2]+...+\alpha[s_r]-\beta[h]=0\ {\rm in}\  H_1(P;{\Z})$$

 Thus there exists a nontrivial class $\eta$ in $H_2(P,\b P;{\Z})$ such that $$\b\eta=((\alpha,\beta),(\alpha,0),...,(\alpha,0),(\alpha,-\beta))$$ 
in $H_1(\b P)=H_1(T_1)\oplus H_1(T_2)\oplus...\oplus H_1(T_{r-1})\oplus H_1(T_r)$. Denote by $H$ an incompressible surface representing $\eta$. Then $H$ is necessarily a horizontal surface and thus $c$ is parallel to some components of $\b H$. This proves the claim.
\end{proof}

\begin{claim}\label{tores}
 The map  $f|T\co T\to N$ is $\pi_1$-injective for any characteristic torus $T$ in $M$. 
\end{claim}
\begin{proof} Let $T$ be a characteristic torus of $M$.  From the Rigidity Theorem of Soma \cite{So} and from hypothesis (i) it is sufficient to consider the case where  $T$ is shared by two distinct Seifert components   
 $\Sigma_1$ and $\Sigma_2$  of $M\setminus\c{T}_M$. For each $i=1,2$, denote by $h_i$ the ${\S}^1$-fiber of $\Sigma_i$. Combining Claims \ref{surface} and \ref{primitive} we deduce that if $f|T\co T\to N$ is not $\pi_1$-injective then $f_{\ast}(h_i)=\{1\}$. Since $h_1$ and $h_2$ generate a rank $2$ subgroup of $\pi_1T$ (by minimality of the JSJ-decomposition) we get a contradiction. This proves the claim. 

\end{proof}
\begin{claim}\label{pull-back}
 There is a map $g$ homotopic to $f$ such that for each Seifert piece of $\Sigma$ of $N$ then each component of $g^{-1}(\Sigma)$ is a Seifert piece of $M$. 
\end{claim}
\begin{proof}
By Lemma \ref{tores} one can apply the  Theorem  \ref{ndeg} combined with \cite[Rigidity Theorem]{So} which imply that one can   arrange $f$ by a homotopy so that for each canonical torus $U$ of $N$ then    $f^{-1}(U)$ is a disjoint union of canonical tori of  $M$. Hence  for each Seifert piece $\Sigma$ of $N$ the space $f^{-1}(\Sigma)$ is a canonical graph submanifold of $M$ (i.e. a submanifold which is the union of some Seifert pieces of $M$).

  If a component $G$ of $f^{-1}(\Sigma)$ is not a Seifert manifold then there exists a canonical torus $T$ in the interior of $G$ which is shared by two distinct Seifert pieces $\Sigma_1$ and $\Sigma_2$ of $G$.  Since by Claim \ref{surface} $f_{\ast}(\pi_1\Sigma_1)$ is not  abelian then using \cite[Addendum to Theorem VI.I.6]{JS} we know that $f|\Sigma_1$ is homotopic to a fiber preserving map. Since $f|T$ is $\pi_1$-injective we get a contradiction by the minimality of the JSJ-decomposiiton.  This proves the claim.
\end{proof}
By hypothesis (ii) one can apply \cite[Rigidity Theorem]{So}.  Thus  one may assume that  $f$ induces a ${\rm deg}(f)$-covering map from $\c{H}(M)$ to $\c{H}(N)$. Moreover since $f$ is $\pi_1$-surjective then to complete the proof of Proposition \ref{reduction} it remains to check the following 
\begin{claim}\label{cov}
There is a map $g$ homotopic to $f$, rel. to $\c{H}(M)$, such that for each Seifert piece $\Sigma$ of $N$ and for each component $G$ of $g^{-1}(\Sigma)$ then $g|G\co G\to\Sigma$ is a covering map.
\end{claim} 
\begin{proof}
First of all note that for each component  $G$ of $f^{-1}(\Sigma)$ then $f|G\co G\to\Sigma$ is non-degenerate and fiber preserving. 

Indeed if $\b G$ is non-empty this comes from Claim \ref{tores}. If $\b G=\emptyset$ then $M=G$ and $N=\Sigma$ so that $f\co M\to N$ is a nonzero degree map between ${\S}^1$-bundles over orientable hyperbolic surfaces. Since ${\rm deg}(f)\not=0$ then $f_{\ast}(h)\not=1$, where $h$ denotes the fiber of $M$ and by Claim \ref{surface}, $f_{\ast}(\pi_1M)$ is a non-abelian group. 

On the other hand, notice that $\Sigma$ is necessarily homeomorphic to a product. 

Indeed if $\b\Sigma\not=\emptyset$ this is obvious and if  $\b\Sigma=\emptyset$ this comes from the following argument: if $\Sigma$ is not homeomorphic to a product then the bundle has a non-zero Euler number and using the Seifert volume in \cite[Theorem 3 and Lemma 3]{Br-G1} and in \cite[Theorem 4]{Br-G2} we get a contradiction (since $G$ has a zero Euler number and ${\rm deg}(f)\not=0)$.

On the other hand notice that for each component  $G$ of $f^{-1}(\Sigma)$ then ${\rm deg}(f|G\co G\to\Sigma)\not=0$.  Indeed, suppose that ${\rm deg}(f|G\co G\to\Sigma)=0$. Since by construction $f|G$ is an allowable map (in the sense of \cite{Ro}) then it induces a zero degree map $\pi\co K\to F$ where $K$, resp. $F$, is a hyperbolic surface such that $G=K\times{\S}^1$, resp. $\Sigma=F\times{\S}^1$. Let $\c{F}$ denote a component of $(f|G)^{-1}(F)$. Arrange $f$ so that $\c{F}$ is incompressible in $G$. Since $f$ is non-degenerate and fiber preserving then $\c{F}$ is necessarily a horizontal surface. Since $f|\c{F}\co\c{F}\to F$ factors throught $\pi$ then ${\rm deg}(f|\c{F}\co\c{F}\to F)=0$. Thus using the naturality property we get the equality 
$$f_{\sharp}(\hat{\Theta}_G[\c{F}])=\hat{\Theta}_{\Sigma}f_{\sharp}[\c{F}]=0$$  
This gives a contradiction with hypothesis (iii) of Proposition \ref{reduction} since $$\|\hat{\Theta}_G[\c{F}]\|_1=\|\c{F}\|>0$$
Hence since  ${\rm deg}(f|G\co G\to\Sigma)\not=0$ then $f$ induces a nonzero degree map $f|\c{F}\co\c{F}\to F$ so that we get 
$$f_{\sharp}(\hat{\Theta}_G[\c{F}])=\hat{\Theta}_{\Sigma}f_{\sharp}[\c{F}]={\rm deg}(f|\c{F}\co\c{F}\to F)\hat{\Theta}_{\Sigma}[F]$$
which implies that 
$$\|\c{F}\|=\left|{\rm deg}(f|\c{F}\co\c{F}\to F)\right|\times\|F\|$$
Thus we get the equality 
$$\|K\|={\rm deg}(\pi)\times\|F\|$$
Hence $\pi$ is homotopic to a covering map which implies that $f|G$ is also homotopic to a covering map. This proves the claim and completes the proof of the Proposition \ref{reduction}.
\end{proof}
\subsection{Proof of Theorem \ref{localiso}}\label{fraks}
First of all note that according to \cite{W}, \cite{So} and the hypothesis of the theorem,  we may assume that either $\c{T}_M\not=\emptyset$ or $M$ is a ${\Hi}^2\times{\R}$-manifold. 
Thus  we have $\hat{H}_2^{l_1}(M;{\R})\not=0$. Then either $N$ has a non-empty JSJ-decomposition or $N$ admits one of the following geometry: ${\Hi}^2\times{\R}$ or $\t{\rm SL}(2,{\R})$.

 Denote by $P_1,...,P_k$ the Seifert pieces of $M$. Then $[M]^{\c{T}_M}=\hat{\Theta}_{P_1}\alpha_{P_1}+...+\hat{\Theta}_{P_k}\alpha_{P_k}$. Then using the additivity property of the $l_1$-norm and the isometry hypothesis we have
$$\sum_{i=1}^k\left\|\hat{\Theta}_{P_i}\alpha_{P_i}\right\|_1\geq\sum_{i=1}^k\left\|f_{\sharp}\hat{\Theta}_{P_i}\alpha_{P_i}\right\|_1\geq\left\|f_{\sharp}\sum_{i=1}^k\hat{\Theta}_{P_i}\alpha_{P_i}\right\|_1=\left\|\sum_{i=1}^k\hat{\Theta}_{P_i}\alpha_{P_i}\right\|_1$$ 
and since, by the additivity property of Theorem \ref{1}, $$\left\|\sum_{i=1}^k\hat{\Theta}_{P_i}\alpha_{P_i}\right\|_1=\sum_{i=1}^k\left\|\hat{\Theta}_{P_i}\alpha_{P_i}\right\|_1$$
hence we get 
$$\sum_{i=1}^k\left\|f_{\sharp}\hat{\Theta}_{P_i}\alpha_{P_i}\right\|_1=\sum_{i=1}^k\left\|\hat{\Theta}_{P_i}\alpha_{P_i}\right\|_1$$
and thus $\|f_{\sharp}\hat{\Theta}_{P_i}\alpha_{P_i}\|_1=\|\hat{\Theta}_{P_i}\alpha_{P_i}\|_1$ for any $i=1,...,k$ which implies, using Lemma \ref{generator}, that 
$\|f_{\sharp}\hat{\Theta}_{P}\alpha\|_1=\|\hat{\Theta}_{P}\alpha\|_1$ for each $\alpha\in\hat{H}_2(P,\b P)$ and $P$ in $\c{S}(M)$. Consider now the following commutative diagram
$$\xymatrix{
M_2 \ar[r]^{f_2} \ar[d]_q & N_2 \ar[d]^p \\
M_1 \ar[r]^{f_1} \ar[d]_s & N_1 \ar[d]^r \\
M \ar[r]^f & N
}$$
obtained as follows. The map $s\co M_1\to M$ is a finite covering such that each Seifert piece of $M_1$ is a product with  at least two boundary components, if $\c{T}_M\not=\emptyset$,  and each canonical torus of $M_1$ is shared by two distinct components of $M_1\setminus\c{T}_{M_1}$ (for the existence of such a covering see \cite[Lemmas 3.2 and 3.5]{DW}), the map  $r\co N_1\to N$ is a finite covering corresponding to the subgroup $f_{\ast}s_{\ast}(\pi_1M_1)$ in $\pi_1N$, which is of finite index since ${\rm deg}(f)\not=0$, the map $f_1\co M_1\to N_1$ is a lifting of $f\circ s$ which exists by our construction, the map $p\co N_2\to N_1$ is a finite covering such that each Seifert piece of $N_2$ is a ${\S}^1$-bundle over an orientable surface and $f_2\co M_2\to N_2$ is the finite covering of $f_1$ corresponding to $p$, and $q\co M_2\to M_1$ is the covering corresponding to the subgroup $(f_1)_{\ast}^{-1}(p_{\ast}\pi_1N_2)$. Notice that it follows from the construction that $f_1$ and $f_2$ are $\pi_1$-surjective. On the other and let $\alpha$ be an element of $\hat{H}_2(P,\b P)$, where $P$ is a Seifert piece of $M_2$. Then using the isometric properties of finite coverings of Proposition \ref{2} together with the commutativity of the diagram we get the following equalities (**)
$$\left\|\hat{\Theta}_P\alpha\right\|_1=\left\|f_{\sharp}s_{\sharp}q_{\sharp}\hat{\Theta}_P\alpha\right\|_1=\left\|r_{\sharp}p_{\sharp}(f_2)_{\sharp}\hat{\Theta}_P\alpha\right\|_1\leq\left\|(f_2)_{\sharp}\hat{\Theta}_P\alpha\right\|_1\leq\left\|\hat{\Theta}_P\alpha\right\|_1$$  
Indeed replacing $\alpha$ by $-\alpha$ we may assume that $\hat{\Theta}_P\alpha\in\hat{H}_2^{l_1,+}\left(M_2^{\c{T}_{M_2}}\right)$ and since the covering maps $q$ and $s$ preserve the torus decompositions then   $q_{\sharp}\hat{\Theta}_P\alpha\in\hat{H}_2^{l_1,+}\left(M_1^{\c{T}_{M_1}}\right)$ and $s_{\sharp}q_{\sharp}\hat{\Theta}_P\alpha\in{\rm Im}(\hat{\Theta}_Q)$ for some Seifert piece $Q$ of $M$. 
Thus one can apply Proposition \ref{reduction} which implies that $f_2$ is homotopic to a homeomorphism. Since $p,q,r,s$ are finite covering maps then this implies that $f$ is $\pi_1$-injective. Consider the finite covering $\t{N}\to N$ corresponding to $f_{\ast}(\pi_1M)$. Then $f$ lifts to a map $\t{f}\co M\to\t{N}$ inducing an isomorphism at the $\pi_1$-level. We deduce from this point using \cite{Wa} and \cite[Section 5.3, Theorem 6]{O} that $\t{f}$ is a homeomorphism. This implies that $f$ is a covering map and completes the proof of Theorem \ref{localiso}.

\subsection{Proof of Corollary \ref{rigidity}} Using the duality we may assume that $f$ induces an isometry $f_{\sharp}\co\hat{H}_2^{l_1}(M;{\R})\to\hat{H}_2^{l_1}(N;{\R})$. Then according to Theorem \ref{localiso} it is sufficient to consider the case 
 where $M$ is a $\t{\rm SL}(2,{\R})$-manifold. Then in particular we have $\hat{H}_2^{l_1}(M;{\R})\not=0$ and thus $\hat{H}_2^{l_1}(N;{\R})\not=0$. On the other hand, since $f\co M\to N$ is a  degree  one map then $N$ admits necessarily a geometry ${\Hi}^2\times{\R}$ or $\t{\rm SL}(2,{\R})$.  Since $f$ is a degree one map, then it is $\pi_1$-surjective and thus it is homotopic to a fiber preserving map and $f_{\ast}([\gamma])\not=\{1\}$ in $\pi_1N$ for any fiber $\gamma$ of $M$. If $N$ is a ${\Hi}^2\times{\R}$-manifold, choose a horizontal surface $H$ is $N$. Then one can arrange $f$ by a homotopy so that $f^{-1}(H)$ is either a horizontal or a vertical surface if $M$. Since $M$ is a $\t{\rm SL}(2,{\R})$-manifold then there are no horizontal surface in $M$ so that $f^{-1}(H)$ consists necessarily of vertical surfaces.  A contradiction since $f$ is a non-degenerate fiber preserving. Thus $N$ is also a $\t{\rm SL}(2,{\R})$-manifold.

Since $N$ is a Haken manifold then there exists a non-empty torus splitting $\c{U}$ of $N$. On the other hand note that one can arrange $f$ by a homotopy so that $\c{T}=f^{-1}(\c{U})$ is also a torus splitting of $M$.  As in paragraph \ref{fraks}, consider  the  commutative diagram
$$\xymatrix{
(M_2,\c{T}_2) \ar[r]^{f_2} \ar[d]_q & (N_2,\c{U}_2) \ar[d]^p \\
(M_1,\c{T}_1) \ar[r]^{f_1} \ar[d]_s & (N_1,\c{U}_1) \ar[d]^r \\
(M,\c{T}) \ar[r]^f & (N,\c{U})
}$$
where  $s\co M_1\to M$ is a finite covering such that each component of $M_1\setminus\c{T}_1$ is a product  with  at least two boundary components, $\c{T}_1=s^{-1}(\c{T})$, and each component of $\c{T}_1$ is shared by two distinct components of $M_1\setminus\c{T}_{1}$, the map  $r\co N_1\to N$ is a finite covering corresponding to the subgroup $f_{\ast}s_{\ast}(\pi_1M_1)$ in $\pi_1N$ and $\c{U}_1=r^{-1}(\c{U})$, the map $f_1\co M_1\to N_1$ is a lifting of $f\circ s$, the map $p\co N_2\to N_1$ is a finite covering with $\c{U}_2=p^{-1}(\c{U}_1)$ such that each component of $N_2\setminus\c{U}_2$ is a product and $f_2\co M_2\to N_2$ is the finite covering of $f_1$ corresponding to $p$ with $\c{T}_2=q^{-1}(\c{T}_1)$. Note that, using the same argument as in (**) paragraph \ref{fraks}, we obtain  $\|(f_2){\sharp}\hat{\Theta}_P({\alpha})\|_1=\|\hat{\Theta}_P\alpha\|_1$ for each $\alpha\in\hat{H}_2(P,\b P)$ when $P$ runs over the  components of $M_2\setminus\c{T}_2$. Hence, since using Claim \ref{surface} and \ref{primitive} the map $f_2|\c{T}_2$ is $\pi_1$-injective,  one can apply the proof of Claim \ref{cov} which implies that for each component $Q$ of $N_2\setminus\c{U}_2$ then $f_2|f_2^{-1}(Q)\co f_2^{-1}(Q)\to Q$ is a covering map. This completes the proof of Corollary \ref{rigidity}.

\end{document}